\documentclass{elsarticle}

\usepackage[utf8]{inputenc}
\usepackage{graphicx} 
\usepackage{caption}
\usepackage{enumitem}

\usepackage{algorithm}
\usepackage[noend]{algpseudocode}
\usepackage{subcaption}
\usepackage{appendix}
\usepackage{anysize}
\usepackage{mathtools}
\usepackage{multicol}
\usepackage{breqn} 
\usepackage{soul}
\usepackage{placeins}
\usepackage{xcolor}
\usepackage{colortbl}
\usepackage{ulem}
\usepackage[rightcaption]{sidecap}

\usepackage{nomencl}
\usepackage{wrapfig}
\usepackage{multicol}
\usepackage{geometry}
\usepackage{comment}
\usepackage{pgfplots}
\usepackage{amsmath}
\usepackage{amssymb}
\usepackage{amsthm}
\usepackage{cancel}
\usepackage{ulem}
\usepackage{graphicx}
\usepackage{booktabs}
\usepackage{xspace} 
\usepackage{nicematrix}
\usepackage{scalerel,stackengine}
\usepackage[colorlinks=true, allcolors=blue]{hyperref}
\usepackage{tikz}

\usetikzlibrary{fit}
\usetikzlibrary{positioning}
\usetikzlibrary{calc}
\usepackage{cleveref}

\newcommand{\coo}{\ensuremath{\mathrm{CO_2}}\xspace}
\newcommand{\Hess}{\mathbb{H}}
\newcommand{\HessNew}{\mathbb{H}}

\newcommand{\OurStabRes}{\text{Current}}
\newcommand{\LitStabRes}{\text{Nichita}}
\newcommand{\OptFunc}{\mathcal{F}}
\newcommand{\doubleEnter}{\\}
\newcommand{\Grad}{g}
\newcommand{\TriangleSide}{1.5}
\newcommand{\x}{\mathbf{x}}
\newcommand{\mycolor}{red}
\newcommand{\concUnit}{mol/m$^3$}
\newcommand{\volUnit}{cm$^3$}
\newcommand{\cons}{\mathcal{C}}
\newcommand{\xOld}{\mathbf{x}}
\newcommand{\xNew}{\mathbf{x}}
\newcommand{\NExpanded}[1]{ N_1^{#1}, \dots, N_n^{#1}}
\newcommand{\NBold}[1]{\mathbf{N}^{#1}}
\newcommand{\NExpandedNoSuper}{ N_1, \dots, N_n}
\newcommand{\NBoldNoSuper}{\mathbf{N}} 
\newcommand{\Res}{(\xi)}
\newcommand{\Smejkal}{\text{UVN}}
\newcommand{\uvn}{\text{UVN\space}}

\newcommand{\tvnArgs}{T, V, \text{ and } \mathbf{N}}
\newcommand{\Ours}{\text{TVN}}

\newcommand{\renderStyle}{\displaystyle}
\newcommand{\myref}{\Cref}
\newcommand{\myeqref}{\Cref}
\newcommand{\NDiff}{\mathbf{N}^{\Res}}

\newcommand{\NStarDiffExpanded}{N_1^\star - \sum_{k=1}^{p-1} N_1^{(k)}, \dots, N_n^\star - \sum_{k=1}^{p-1} N_n^{(k)}}
\newcommand{\vertVec}[3]{\begin{array}{c} \renderStyle  #1 #3 \renderStyle  \vdots #3 \renderStyle  #2 \end{array}}

\newcommand{\pd}[2]{\frac{\partial #1}{\partial #2}}
\newcommand{\pdd}[3]{\frac{\partial^2 #1}{\partial #2 \partial #3}}

\newcommand{\colvec}[1]{\begin{pmatrix} #1 \end{pmatrix}}
\newcommand{\Lag}{\mathcal{L}}
\newcommand{\SCL}{S\cons\Lag}
\newcommand{\ACL}{A\cons\Lag}

\allowdisplaybreaks
\stackMath

\title{A Fast and Robust Reformulation of the UVN-Flash Problem via Direct Entropy Maximization}

\makenomenclature
\begin{document}

\begin{frontmatter}


\author[delft,cwi]{Pardeep Kumar\corref{cor1}}
\ead{pardeep@cwi.nl}

\author[shell]{Patricio I. Rosen Esquivel}

\cortext[cor1]{Corresponding author}

\affiliation[delft]{
    organization={Delft University of Technology},
    city={Delft},
    country={The Netherlands}
}

\affiliation[cwi]{
    organization={Centrum Wiskunde \& Informatica},
    city={Amsterdam},
    country={The Netherlands}
}

\affiliation[shell]{
    organization={Shell Projects and Technology},
    city={Amsterdam},
    country={The Netherlands}
}

\begin{abstract}
We investigate the phase equilibrium problem for multicomponent mixtures under specified internal energy (\( U \)), volume (\( V \)), and mole numbers (\( N_1, N_2, \dots, N_n \)), commonly known as the UVN-flash problem. While conventional phase equilibrium calculations typically use pressure-temperature-mole number (\( PTN \)) specifications, the UVN formulation is essential for dynamic simulations of closed systems and energy balance computations. Existing approaches, including those based on iterative pressure-temperature updates and direct entropy maximization, suffer from computational inefficiencies due to nested iterations and reliance on inner Newton solvers.  

In this work, we present a novel reformulation of the UVN-flash problem as a direct entropy maximization problem that eliminates the need for inner Newton iterations, addressing key computational bottlenecks. We derive two new novel formulations: 1) a formulation based on entropy and internal energy and (2) an alternative formulation based on Helmholtz free energy. We begin with a stability analysis framework, followed by a reformulation of the UVN flash problem in natural variables. We then introduce our novel approach and discuss the numerical methods used, including gradient and Hessian computations. The proposed method is validated against benchmark cases, demonstrating improved efficiency and robustness.

\end{abstract}




\end{frontmatter}


\section{Introduction}

We investigate the phase equilibrium calculations for multicomponent mixtures under specified internal energy ($U$), volume ($V$), and mole numbers ($N_1, N_2, \dots, N_n$), commonly referred to in the literature as the UVN stability test or the UVN-flash problem. Compared to the more conventional PTN flash (e.g., \cite{michelsen_isothermal_1981, michelsen_isothermal_1982, michelsen_thermodynamic_2007}) (where pressure, temperature, and mole numbers are specified), the UVN specification is less commonly addressed. However, it plays a crucial role in various thermodynamic applications where energy and volume are conserved, such as in the dynamic simulation of closed systems and energy balance calculations in process design. Notably, the UVN-flash formulation proves to be particularly valuable in non-isothermal problems, such as those encountered in the dynamic simulation of tanks and \coo injection in geological storage. Key contributions in this area include~\cite{arendsen_dynamic_2009, castier_dynamic_2010,goncalves_dynamic_2007,lu_transient_2014, muller_dynamic_1997}. 

Additionally, the UVN-flash problem has also been explored in the work of some researchers, albeit not as the primary focus. For instance, Lipovac et al.~\cite{lipovac_unified_2024} developed a unified flash procedure for isenthalpic (PHN) and isochoric (UVN) flash calculations, where a persistent set of unknowns and equations is solved during equilibrium calculations, enabling simultaneous phase stability and split calculations. Fathi~\cite{fathi_rapid_2021} examined flash calculations using volume functions and reduction methods, which can be extended to UVN scenarios. 

Michelsen~\cite{michelsen_state_1999} proposed a general framework to address flash problems under various specifications, including UVN. His approach utilizes the PTN-flash in an inner loop while iteratively updating pressure and temperature in an outer loop. The advantage of this method is that it leverages existing PTN-flash solvers, but the nested iterations quickly render it computationally expensive.

One of the earliest works addressing the UVN-flash problem is by Saha et al.~\cite{saha_isoenergetic-isochoric_1997}. In this paper, the authors develop heuristics to estimate pressure and temperature corresponding to specified UVN conditions. Their approach combines successive substitution (fixed-point iteration) for updating equilibrium $K$-values with Newton's method for pressure and temperature updates. However, they often encountered convergence to trivial solutions, limiting the robustness of their method.

Bi et al.~\cite{bi_efficient_2020} reformulate the UVN-flash problem using the Rachford-Rice equation while ensuring pressure equilibrium and enforcing internal energy and volume constraints. Their approach employs fixed-point iteration with soft tolerance, followed by Newton’s method for refinement. However, Nichita~\cite{nichita_robustness_2023} has highlighted certain limitations in its robustness.

A significant contribution on UVN flash is found in the work of Castier~\cite{castier_solution_2009}, who proposed direct entropy maximization as an alternative approach. In his method, the algorithm adaptively adds or removes phases as needed during the computation. However, obtaining a good initial phase split requires a reasonable estimate of pressure and temperature, which has to be determined using heuristics. In cases of numerical difficulties, Castier's method switches to a PTN-flash solver for the inner loop while adjusting pressure and temperature in the outer loop, ensuring that $U$ and $V$ approach their specified values. Once sufficient estimates for $P$ and $T$ are found, the algorithm returns to direct entropy maximization.

Another important contribution is by Smejkal et al~\cite{smejkal_phase_2017}, who also applied direct entropy maximization for both stability and flash calculations. They used the results of stability analysis as initial guesses for the flash calculations, demonstrating the utility of entropy-based methods in UVN-flash scenarios. This approach requires an inner Newton iteration to determine temperature by solving \(U(T, V, \NBoldNoSuper) = U\) for given \(U, V, \NBoldNoSuper\) which poses an additional computational burden.

In this work, we revisit the UVN-flash problem as direct entropy maximization problem and propose a novel reformulation that eliminates the need for inner Newton iterations, thus removing the computational bottleneck encountered in the traditional UVN approach via entropy maximization. The structure of the paper is as follows. We begin with a precursor to stability analysis in  \myref{sec:stability_analysis}, where we present the relevant formulations and the generation of initial guesses for both stability and flash calculations. This is followed by a recap of the UVN flash formulation in natural variables in \myref{sec:Mikyska_UVN}. Next, in \myref{sec:novel_contribution}, we introduce our novel contributions by reformulating the UVN-flash problem. Building on this framework, we develop two novel formulations: (1) a formulation based on entropy and internal energy and (2) an alternative formulation based on Helmholtz free energy in \Cref{sec:numerical_method}. We also discuss the numerical approach in \myref{sec:numerical_method}, including the necessary gradient and Hessian computations (derived in \ref{sec:gradAndHessNew}). Finally, we present the results in \myref{sec:results} and conclude with key findings and implications in \myref{sec:conclusion}.

\section{Preliminaries}

For sake of clarity, in this section, we define the following concepts in the context of UVN-flash.
\subsection{Trial Phase:} The trial phase is an incipient phase introduced to assess the thermodynamic stability of a system. It involves perturbing the composition of the system and evaluating whether the introduction of this new phase leads to an increase in entropy (for UVN flash calculations). If the entropy increases, the system is unstable as a single phase, and phase separation is favorable. 

\subsection{Reference Phase:}
The reference phase (\(\star\)) represents a hypothetical single-phase system characterized by the total internal energy \( U^{\star} \), volume \( V^{\star} \) and total mole numbers \( \NBold{\star} = (\NExpanded{\star}) \). Introducing a trial phase forms a two-phase system. Stability is assessed by comparing the entropy of the two-phase system with the reference phase. If the two-phase system has more entropy than the reference phase, then the reference phase is deemed unstable, and phase separation occurs.

\subsection{Stability Analysis:} 
Stability analysis is the first step in flash calculations, as it determines the stability of a multicomponent mixture across \( p \) phases, where \(p\) represents the number of phases in the system. The primary objective of this analysis is to establish whether the mixture will remain stable as a single phase or if it will separate into \( p + 1 \) phases. In most cases, we focus on systems with at most two phases, typically a vapor-liquid mixture. In such cases, stability analysis determines whether the mixture can remain as a single phase or will separate into a vapor-liquid equilibrium. 

A crucial aspect of stability analysis is its role in providing an initial guess for subsequent flash calculations. If instability is detected, the analysis often yields valuable information about the incipient phase, such as its temperature, concentration, and internal energy density, which can significantly aid in the convergence of the flash calculation algorithm. 

\subsection{Flash Calculation}  
When a stability test indicates that a mixture is thermodynamically unstable, a flash calculation is performed to determine the phase equilibrium of the multicomponent mixture under specified conditions, such as pressure and temperature, internal energy and volume, or entropy and volume. The flash calculation predicts the amounts and compositions of each phase, assuming the system reaches equilibrium.  

Consider a mixture of \( n \) components at a given internal energy \(U\), volume \( V \) and the composition \( \NBoldNoSuper = \{ \NExpandedNoSuper \} \), where \( N_i \) is the mole number of the \( i^{\text{th}} \) component. Assuming the mixture exists in vapor-liquid equilibrium, the objective of the flash problem is to determine the phase split, i.e., the number of moles of each component in the vapor and liquid phases, along with the thermodynamic properties (such as the temperature, internal energy and enthalpy etc.) of each phase. Flash calculations are also referred to as phase split or phase equilibrium calculations.

In this study, the UVN flash problem is addressed using an equation of state derived from the Helmholtz energy function. Specifically, we will use Peng-Robinson \cite{peng_new_1976} Equation of State (EOS) for all our results. This EOS enables the computation of thermodynamic properties for given \( T \), \( V \), and \( \NBoldNoSuper \). For specified values of internal energy \( U \), volume \( V \), and mole numbers \( \NBoldNoSuper \), the phase split calculations proceed as follows: 

\begin{enumerate}  
    \item Perform a stability analysis to determine whether the mixture is stable or if phase separation occurs. If the mixture is unstable, the results from the stability analysis provide an initial guess for the flash calculations.  
    \item If the mixture is unstable, initiate flash calculations to determine the equilibrium temperature \( T \), the phase volumes \( V^{(1)}, V^{(2)} \dots V^{(p)} \) (where the superscripts refer to the phase index), and the mole numbers of each component in each phase.  
\end{enumerate}  
Accordingly, we begin with stability analysis in the next section, followed by the flash procedure in the subsequent section.

\section{Stability Analysis} \label{sec:stability_analysis}

Stability analysis is a fundamental step in assessing the thermodynamic stability of a mixture and determining whether phase separation occurs. In this section, we present the formulation of the UVN stability problem, discuss initialization strategies, and outline an algorithm to generate initial guesses for flash calculations based on stability analysis results.

The UVN stability problem can be reduced to the VTN stability problem as follows. For given \(U^\star, V^\star\)  and \(\mathbf{N}^{\star}\), we can solve
\begin{equation}
    U(T, V^{\star}, \mathbf{N}^{\star}) = U^{\star},
    \label{eqn:internal_energy_eq}
\end{equation}  
for \(T\), as discussed by Mikyska~\cite{mikyska_investigation_2012} and Nichita~\cite{nichita_unified_2024}. Therefore, we provide the formulation of VTN stability analysis in the following subsection.

\subsection{VTN Stability Formulation}  
In this section, we briefly discuss the formulation of VTN stability as UVN stability can be reduced to VTN stability~\cite{nichita_robustness_2023}. For a detailed derivation of the VTN stability formulation, we refer to the work of Mikyska et al.~\cite{mikyska_investigation_2012}. Mikyska introduced the volume function \( \Phi \) to express the VTN stability condition, which is formulated as follows:  
\begin{align}
    \ln \frac{c^\prime x_i}{cz_i} + \ln \Phi_i(\boldsymbol{cz}) - \ln \Phi_i(\boldsymbol{c^\prime x}) = 0, \quad \forall i \in \{1, \dots, n\}
    \label{eqn:stability_vol_func}
\end{align}
where \(n\) is the number of components, \(c^\prime = \frac{N^\prime}{V^\prime}\) is the molar concentration in the trial phase, \(z_i = \frac{N_i}{N}\) and \(x_i = \frac{N^\prime_i}{N^\prime}\) are the mole fractions of component \(i\) in the reference and trial phases, respectively. Here, \(\boldsymbol{cz} = \{cz_1, cz_2, \dots, cz_n\}\) and \(\boldsymbol{c^\prime x} = \{c^\prime x_1, c^\prime x_2, \dots, c^\prime x_n\}\) are the molar concentration vector in the reference phase and the the trial phase respectively. The volume function of the component \(i\), \(\Phi_i\) is related to the fugacity coefficient of the component \(i\), \(\phi_i\) via 

\begin{equation}
    \Phi_i = \frac{1}{Z \phi_i},    
\end{equation}
where \(Z = \frac{PV}{nRT}\) is the compressibility factor, \(R\) is the universal gas constant, \(P\) is the pressure, \(T\) is the temperature, and \(V\) is the volume.

To perform the stability analysis, the system of \(n\) equations in \eqref{eqn:stability_vol_func} must be solved for the trial-phase molar concentrations, \( \boldsymbol{c^\prime x} \). These concentrations are used to evaluate the tangent plane distance (TPD), denoted by \( D \), as given by:

\begin{equation} \label{eqn:TPD_func}
    D = u^\prime\left(\frac{1}{T^\prime} - \frac{1}{T^\star}\right) + \left(\frac{P^\prime}{T^\prime} - \frac{P^\star}{T^\star}\right) - \sum_{i=1}^n \left(\frac{\mu^\prime_i}{T^\prime} - \frac{\mu^\star_i}{T^\star}\right)c_i^\prime.
\end{equation}
where \( u^\prime \) is the internal energy density of the trial phase,  \( T^\prime \) and \( T^\star \) are the temperatures in the trial and reference phases, respectively, \( P^\prime \) and \( P^\star \) are the pressures in the trial and reference phases, respectively, \( \mu^\prime_i \) and \( \mu^\star_i \) represent the chemical potentials of component \( i \) in the trial phase and the reference phase, respectively and \( c_i^\prime \) is the molar concentration of component \( i \) in the trial phase. A positive value of \( D \) indicates that the initial phase is thermodynamically unstable, implying that phase splitting will occur (see~\cite{smejkal_phase_2017} for more details). The convergence of the stability test strongly depends on the choice of an appropriate initial guess for \( c^\prime x \). A judicious choice of initialization is crucial for ensuring numerical stability and enhancing the robustness of the analysis. We continue in the next section with the initialization strategy.

\subsection{Initialization For Stability Analysis} \label{subsec:simplex_IC} 

In this section, we discuss the initialization strategy for VTN stability analysis. We adopt the simplex-based initialization method proposed by Smejkal et al.~\cite{smejkal_phase_2017}, which leverages the geometric properties of the feasibility domain of admissible molar concentrations. In this approach, the feasible domain is represented as an \(n\)-simplex, where \(n\) denotes the number of components in the mixture. Initial guesses are generated by computing the barycenter of the simplex and the midpoints between the barycenter and each of the \(n + 1\) vertices. This procedure yields \(n + 2\) initial estimates, ensuring a well-distributed set of starting points for the stability analysis.  

\begin{figure}[H]
    \centering
    \begin{tikzpicture}[scale=4]
        \coordinate (V0) at (0,0);
        \coordinate (V1) at (\TriangleSide,0);
        \coordinate (V2) at (0,\TriangleSide);
        
        \path (V0) -- (V1) -- (V2) coordinate (Centroid) at ($1/3*(V0) + 1/3*(V1) + 1/3*(V2)$);
        
        \coordinate (M0) at ($(Centroid)!0.5!(V0)$);
        \coordinate (M1) at ($(Centroid)!0.5!(V1)$);
        \coordinate (M2) at ($(Centroid)!0.5!(V2)$);
        
        \draw [thick, line width=0.75mm] (V0) -- (V1) -- (V2) -- cycle;
        
        \draw [thick, dash pattern=on 6pt off 3pt, line width=0.75mm] (Centroid) -- (V0);
        \draw [thick, dash pattern=on 6pt off 3pt, line width=0.75mm] (Centroid) -- (V1);
        \draw [thick, dash pattern=on 6pt off 3pt, line width=0.75mm] (Centroid) -- (V2);
        
        \foreach \p in {Centroid, M0,M1,M2} {
            \draw[line width=0.7mm, color=\mycolor] (\p) --++ (0.05,0) --++ (-0.1,0); 
            \draw[line width=0.7mm, color=\mycolor] (\p) --++ (0,0.05) --++ (0,-0.1); 

        }
        
        \node[below left] at (V0) {\large $V_0 \ (0,0)$};
        \node[below right] at (V1) {\large $V_1 \ (1/b_1, 0)$};
        \node[above] at (V2) {\large $V_2 \ (0,1/b_2)$};
        
        \node[below right, xshift=3pt] at (M2) {\small M2};
        \node[below, xshift=5pt, yshift=-3pt] at (M1) {\small M1};
        \node[below, xshift=3pt, yshift=-3pt] at (M0) {\small M0};
        \node[left, xshift=20pt, yshift=2pt] at (Centroid) {\small C};

    \end{tikzpicture}
    \caption{Depiction of the initial guesses for $\mathbf{c'}$ in a binary mixture where the two components have molar volumes $b_1$ and $b_2$, respectively.(Adapted from~\cite{bi_efficient_2020})}
    \label{fig:simplex}
\end{figure}
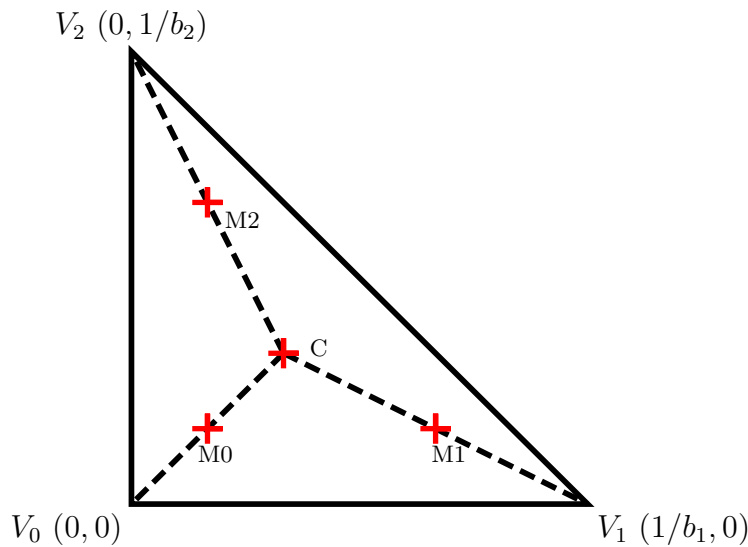

The admissible molar concentrations  $c^{\prime}_{i}$ must satisfy the following conditions:

\begin{align} \label{eqn:feasiblity_cond} 
    \sum_{i=1}^n c^{\prime}_{i}b_i < 1, \quad c^{\prime}_{i} \geq 0, \quad b_i > 0, \quad \forall i \in \{1, \dots, n\},
\end{align}
where  $b_i$ denotes the co-volume of the component $i$ from the Peng-Robinson EOS. Figure \ref{fig:simplex} (adapted from~\cite{bi_efficient_2020}) illustrates the initial concentration guesses for a binary mixture. These initial guesses correspond to four distinct points marked with circles: the barycenter \(C\) and the midpoints \(M_0, M_1,\) and \(M_2\). 
%
These initial guesses serve as starting points for the stability analysis. The results of the stability analysis then are used to generate the initial guess for phase split calculations. These results, however, are in the form of concentration and temperature of the trial phase. A procedure is needed to convert the stability analysis results into the initial guess for phase split calculations, which is addressed in the following subsection.

\subsection{Initial Guess for Flash from Stability Analysis} \label{sec:initial_guess_from_stability}

The results of the stability analysis provide the initial guess required for phase split calculations. A good initial guess is crucial for ensuring convergence in the numerical optimization procedures used in phase split calculations, as discussed in Section \ref{sec:numerical_method}. However, the results from stability analysis are not immediately suitable as initial guesses for flash calculations. Stability analysis provides the concentrations of the incipient phase along with the specific internal energy; but an additional parameter - the volume of the trial phase is needed to initiate phase split calculations. 

We begin by assuming that the trial phase occupies half of the total system volume. The mole numbers of each phase are determined by multiplying the phase volume with the species concentrations obtained from the stability test. The internal energy of each phase is then computed using the internal energy density and phase volume. Next, the phase temperature is determined by solving \Cref{eqn:internal_energy_eq} and finding a temperature consistent with the given internal energy, volume, and mole numbers. This provides a complete initial estimate.

This initial estimate is then iteratively refined by maximizing entropy while simultaneously satisfying the feasibility conditions \eqref{eqn:feasiblity_cond}. At each step, the total entropy of the two-phase system is evaluated. If the entropy increases and all feasibility conditions are met, the solution is accepted. If these criteria are not satisfied, the trial phase volume is further halved, and the internal energy and mole numbers are adjusted accordingly. This iterative process continues until a feasible phase split is achieved or until the predefined iteration limit is reached. This algorithm is outlined in \Cref{alg:init_guess}. 

The final feasible solution obtained from the stability analysis serves as the initial guess for phase split calculations. In the following section, we review the phase split calculation method presented by Castier~\cite{castier_solution_2009} and Smejkal et al.~\cite{smejkal_phase_2017}, which serves as the foundational framework for our work.

\begin{algorithm}[H]
    \caption{Initial Guess Generation and Feasibility Check for Phase Equilibrium}
    \label{alg:init_guess}
    \begin{algorithmic}[1]
    \Require Total internal energy \( U^{\star} \) [J], volume \( V^{\star} \) [m$^3$], and mole numbers \(\mathbf{N}^{\star}=[\NExpanded{\star}] \),  
the trial phase concentration vector \( \mathbf{c} \) [mol/m$^3$]  
and the trial phase internal energy density \( u \) [J/m$^3$]. 
    \Ensure Feasible initial guess for phase split or termination if no solution exists

    \State Compute temperature \( T^\star = \text{for given } U^{\star}, V^{\star}, \mathbf{N}^{\star} \)  of the reference phase by solving \myeqref{eqn:internal_energy_eq}.
    \State Compute entropy \( S^\star = S(T^\star, V^{\star}, \mathbf{N}^{\star}) \) of the reference phase using the equation of state (EOS)
   
   \State Initialize the trial phase \(I\) as follows: 
        \begin{align*}
            V^{I} &= 0.5 \cdot V^{\star} \\ 
            \mathbf{N}^{I} &= V^{I} \cdot \mathbf{c} \\
            U^{I} &= u \cdot V^{I} \\
        \end{align*}
    \State Initialize iteration count: \( n_{\text{iters}} \gets 0 \)

    \While{\(n_{\text{iters}} < \text{maxiters}\)}
        \State Compute total entropy for the two-phase system:
        \[
        S_{\text{two-phase}} = S(U^{I}, V^{I}, \mathbf{N}^{I}) + S(U^{\star} - U^{I}, V^{\star} - V^{I}, \mathbf{N}^{\star} - \mathbf{N}^{I})
        \]
        \State Compute entropy difference: \( \Delta S = S_{\text{two-phase}} - S^\star \)
        \State Update phase properties vector \( \mathbf{x} \) for the trial phase:
        \[
        \mathbf{x} = \begin{bmatrix}
            \mathbf{N}^{I}, & V^{I}, & U^{I}
        \end{bmatrix}
        \]
        \State Check feasibility of \( \mathbf{x} \) using equations \eqref{eqn:feasiblity_cond}.

        \If{\(\Delta S > 0\) \textbf{and} \(\mathbf{x} \text{ is feasible} \)}
            \State \textbf{Return} feasible initial guess \( \mathbf{x} \)
        \EndIf

        \If{\(V^{I} / V^{\star} < 10^{-8}\)}
            \State \textbf{Terminate}: No feasible solution found
        \EndIf

         \State Update phase properties:
            \begin{align*}
                V^{I} &\gets V^{I} / 2 \\
                U^{I} &\gets u \cdot V^{I} \\
                \mathbf{N}^{I} &\gets V^{I} \cdot \mathbf{c} 
            \end{align*}
        \State Increment iteration count: \( n_{\text{iters}} \gets n_{\text{iters}} + 1 \)
    \EndWhile

    \State \textbf{Return} failure: No feasible solution found
    \end{algorithmic}
\end{algorithm}

\section{Direct Entropy Maximization Formulation for UVN Flash Calculations} \label{sec:Mikyska_UVN}


The UVN flash problem can be formulated as a direct entropy maximization problem, constrained by specified system properties as discussed by Castier~\cite{castier_solution_2009} and Smejkal et al.~\cite{smejkal_phase_2017}. Consider a multicomponent mixture composed of $n$ species, distributed across $p$ phases and a total energy $U^\star$, total volume $V^\star$ and the mole numbers vector \(\NBold{\star}\). The total entropy of the system, denoted as $S^{(\Smejkal)}$, can be expressed as:
\begin{align}
    S^{(\Smejkal)} = \sum_{k=1}^{p} S(U^k, V^k, \NBold{k}),
\end{align}
where $U^k$, $V^k$, and $\NBold{k} = \{N_1^{(k)}, \dots, N_n^{(k)}\}$ represent the internal energy, volume, and mole numbers of each component in phase $k$, respectively. The superscript \((\Smejkal)\) highlights the fact that entropy here is expressed as a function of \(U, V, \mathbf{N}\). Additionally, the problem is subject to the following constraints:

\begin{align}
    U^\star = \sum_{k=1}^{p} U^k, \quad V^\star = \sum_{k=1}^{p} V^k, \quad N_i^\star = \sum_{k=1}^{p} N_i^{(k)}, \quad i = 1, \dots, n.
\end{align}
To simplify the problem, we can apply these constraints and reformulate the problem as an unconstrained optimization problem. This is done by writing the properties of phase $p$ as a function of the properties in the other phases. For the entropy function, this reads: 

\begin{align}
    S_{\text{red}}^{(\Smejkal)} = \left[\sum_{k=1}^{p-1} S(U^k, V^k, \NBold{k})\right]  + S(\x^\Res), \label{eqn:entropy_red}
\end{align}
where 
\begin{align} \label{xstar}
\x^\Res := \left(U^\star - \sum_{k=1}^{p-1} U^k, V^\star - \sum_{k=1}^{p-1} V^k, \NDiff \right),  
\end{align}
where \(\NDiff := \{\NStarDiffExpanded\}\). Throughout this text, the subscript \text{red} is used to denote the reduced quantity and the superscript 
\(\Res\) is used to denote the remaining phase, specifically referring to phase 
$p$.

The unconstrained optimization problem now involves solving for the $(p-1)(n+2)$ unknowns: $U^k$, $V^k$, and $\NBold{k}$ for each phase $k \in \{1, 2, \dots, p-1\}$.
Formally, we want to solve the following unconstrained optimization problem:
\begin{align} \label{eqn:opt_prob}
    \x &= \arg\max_{\mathbf{y}} S_{\text{red}}^{(\Smejkal)}({\mathbf{y}}),  
\end{align}
where $\mathbf{y}$ now entails all the $(p-1)(n+2)$ unknowns. The solution to this optimization problem will be discussed in detail in \myref{sec:numerical_method}. We will refer to this approach as the \textbf{\uvn approach} throughout the rest of the paper. 

Before proceeding further, we highlight a difficulty inherent to this approach. Given that the equation of state is in the form \( f(T, V, \NBoldNoSuper) \), \Cref{eqn:entropy_red} requires writing the entropy \(S\) as function of $U$, $V$ and $\NBoldNoSuper$. Therefore, the approach requires determining the temperature \( T \) by solving the equation: 
\[ U(T, V, \NBoldNoSuper) = U \] 
for given \(U, V \text{and} \  \NBoldNoSuper\). Once $T$ is determined, the equation of state can then be used to compute the corresponding entropy. This process results in a nested Newton's method where each iteration of the outer optimization problem requires multiple iterations of the inner solver to achieve convergence, thereby increasing the computational complexity of the solution procedure. We refer to this process as \textbf{EOS inversion} throughout the paper.

To circumvent this difficulty, we reformulate the optimization problem directly in terms of \( \tvnArgs \). This reformulation avoids the need to perform multiple EOS inversions and will be discussed in the next section.

\section{Reformulation of Direct Entropy Maximization: Transition from Unconstrained UVN to Constrained VTN Space} \label{sec:novel_contribution}
This section presents our novel contribution, which relies on the strategy of reformulating the optimization problem in terms of the variables inherent to the Helmholtz energy-based equation of state (EOS), specifically in the \(TVN\)-space. This reformulation circumvents the need to repeatedly invert the EOS to determine the temperature at each iteration. 

The objective function in this new formulation is given by:

\begin{align} \label{eqn:entropy_red_novel}  
    S_{\text{red}}^{(\Ours)} &= \left[ \sum_{k=1}^{p-1} S\left(T, V^{(k)}, \NBold{k}\right)\right]  + S\left(T, V^\star - \sum_{k=1}^{p-1} V^k, \NDiff \right),    
\end{align}
subject to the constraint \( U^\star = \sum_{k=1}^{p} U\left(T, V^{(k)}, \NBold{k}\right)\).  This approach allows for the computation of entropy without the need to explicitly invert the relation \( U(T, V, \NBoldNoSuper) = U \) to determine the temperature \( T \). Note that the constraints for the total volume $V^\star$ and total moles $\NBold{\star}$ are directly incorporated in the arguments of the entropy of phase $p$. However, the constraint of internal energy has not yet been incorporated. Rewriting this constraint in functional form yields

\begin{equation} \label{eqn:cons_red_novel}  
 \cons(\x) := \left[\sum_{k=1}^{p-1} U\left(T, V^{(k)}, \NBold{k}\right) \right] 
    + U\left(T, V^\star - \sum_{k=1}^{p-1} V^k, \NDiff \right) -  U^\star
\end{equation}
The solution \(\hat{\x}\) satisfies the following constrained optimization problem:
\begin{align} \label{eqn:opt_prob_new} \max_{\x} \quad S_{\text{red}}^{(\Ours)}(\x) \quad \text{subject to } \cons (\x) = 0. \end{align}
where the optimization variable $\x := \left(T, V^{(1)}, \NBold{1} , \dots, V^{(p-2)}, \NBold{p-2}, \dots, V^{(p-1)}, \NBold{p-1} \right)$ is a vector of $(p-1)(n+1) + 1$ unknowns.
This constrained optimization problem can be reformulated as an unconstrained one using the method of Lagrange multipliers (see \ref{app:lagrange_mul}). The Lagrangian function is defined as:
\begin{align} \label{eqn:lagrangian}
\mathcal{L}(\x, \lambda) = S_{\text{red}}(\x) + \lambda \, \cons(\x)
\end{align}
where \(\lambda\) is the Lagrange multiplier. For convenience, where possible we omit the superscript \( (\Ours) \) from the objective function. However, it is included when needed to ensure clarity. To find the optimum, we solve the following system of equations:

\[
\nabla_{\x, \lambda} \mathcal{L} = 0
\]
where \( \nabla_{\x, \lambda} \mathcal{L} = \left( \nabla_{\x} \mathcal{L}, \frac{\partial \mathcal{L}}{\partial \lambda} \right) \). The gradient $\nabla_{\x} \mathcal{L}$ of the Lagrangian with respect to \( \x \) is given by:

\[
\nabla_{\x} \mathcal{L} = \left( \frac{\partial \mathcal{L}}{\partial \x_1}, \frac{\partial \mathcal{L}}{\partial \x_2}, \dots, \frac{\partial \mathcal{L}}{\partial \x_{(p-1)(n+1) + 1}} \right)
\]
The condition \( \nabla_{\x, \lambda} \mathcal{L} = 0 \) leads to two sets of equations.

\begin{enumerate}
    \item Stationarity Condition: 
        \begin{equation} \label{eqn:kkt1}
            \nabla {S_{\text{red}}}({\x}) = -\lambda \nabla \cons(\x)
        \end{equation} 
        This ensures that the gradient of the objective function \(S_{\text{red}}\) is parallel to the gradient of the constraint \(\cons\).

    \item Primal Feasibility Condition:
        \begin{equation} \label{eqn:kkt2}
            \cons(\x) = 0
        \end{equation}
        This ensures that the constraint is satisfied.
\end{enumerate}
This formulation leads to a system of $(p-1)(n+1) + 2$ equations. Specifically, for $p=2$, the system contains one additional equation compared to the approach of Smejkal et al. \cite{smejkal_phase_2017}. For $p=3$, the number of equations is the same in both approaches. However, for ( $p \geq 4$ ), our approach requires solving ($p-3$) fewer equations as compared to Smejkal et al. \cite{smejkal_phase_2017}. 

In summary, this section has outlined the framework of our novel approach. The next section begins with a discussion of the numerical optimization procedure, followed by a simplification of the Lagrangian and the introduction of two new formulations.

\section{Numerical Optimization} \label{sec:numerical_method}
This section outlines the optimization of the objective function, as defined in \myeqref{eqn:entropy_red} for Smejkal's approach and \myeqref{eqn:lagrangian} for our method. Additionally, we derive an explicit expression for the Lagrange multiplier and use it to simplify the Lagrangian in our approach, making the implementation more straightforward. Finally, we conclude this section by performing a consistency check of our new formulation. For simplicity, we assume the number of phases is known a priori, as determined by a stability analysis. 
\subsection{Optimization of the objective function} \label{sec:newton_method}
Formally, we seek to solve the following unconstrained optimization problem:

\begin{equation}
    \xNew = \arg\max_{\mathbf{y}} \OptFunc({\mathbf{y}})
\end{equation}
where $\OptFunc$ is the objective function. To solve this, we need to find the critical points of the gradient of the objective function, denoted as \( \Grad(\xOld) \). This gradient is expressed as:

\begin{align} \label{eqn:uncons_entropy_grad}
    \Grad(\x) = \begin{cases}
       \renderStyle \nabla S_{\text{red}}^{(\Smejkal)}({\x}),  \quad & \uvn \\ \\
        \renderStyle \nabla S_{\text{red}}^{(\Ours)}({\x}) + \lambda \ \nabla \cons(\x), \quad & \text{Ours}
    \end{cases}     
\end{align}
We will revisit the alternate forms of \(\Grad(\x)\) for our approach in the next section where we derive the Lagrange Multiplier \(\lambda\). The specific forms of the optimizer $\x$ are given by:

\begin{align} \label{eqn:optimisers}
    \x = \begin{cases}
         (\NBold{1}, V^{(1)}, U^{(1)}, \dots, \NBold{p-1}, V^{(p-1)}, U^{(p-1)}), \quad &\uvn \\ \\
         (\NBold{1}, V^{(1)}, \dots, \NBold{p-2}, V^{(p-2)}, \NBold{p-1}, V^{(p-1)}, T), \quad &\text{Ours} 
    \end{cases}     
\end{align}
Here, \( N_i^{(k)} \) represents the mole number of component \( i \) in phase \( k \), while \( U^{(k)} \) and \( V^{(k)} \) correspond to the internal energy and volume of phase \( k \), respectively.

The optimization problem can now be written as solving \( \Grad(\xOld) = 0 \). This is a nonlinear system, and it can be solved using a nonlinear solver, such as Newton-Raphson or a variant. To apply the Newton-Raphson method, we need the gradient of \( \Grad(\xOld) \). The gradient  of \( \Grad(\xOld) \) is the Hessian of the objective function \(\OptFunc\), given by:

\[
\Hess(\x) = \left[ \pdd{\OptFunc}{\x_i}{\x_j} \right].
\]
In this context, the Hessian for both optimization approaches is expressed as:

\begin{align}
   \Hess(\x) = \begin{cases}
        \displaystyle \pdd{S_{\text{red}}^{(\Smejkal)}(\x)}{\x_i}{\x_j}, & \quad \uvn \\\\
        \displaystyle \pdd{S_{\text{red}}^{(\Ours)}(\x)}{\x_i}{\x_j} + \lambda \pdd{\cons(\x)}{\x_i}{\x_j}, & \quad \text{Ours} 
    \end{cases}
\end{align}
where \( \x \) is defined as per \myeqref{eqn:optimisers}, and the entropy function \( S_{\text{red}}(\x) \) is given by \myeqref{eqn:entropy_red} in Smejkal's formulation and by \myeqref{eqn:entropy_red_novel} in our approach. The additional term in our formulation accounts for the contribution of the constraint function \( \cons(\x) \) through the Lagrange multiplier \( \lambda \), ensuring that the optimization respects the imposed constraints.

\myeqref{eqn:uncons_entropy_grad} can be solved using a non-linear solver. We employ Newton's method, which updates the solution iteratively as follows:

\begin{equation} \label{eqn:newton_raphson}
    \mathbf{x}_{k+1} = \mathbf{x}_k + \alpha_k \Delta \mathbf{x}_k    
\end{equation}
where \( \alpha_k \) is the step size and the update direction \( \Delta \mathbf{x}_k \) satisfies:

\begin{equation} \label{eqn:linear_system}
    \Hess(\mathbf{x}_k) \Delta \mathbf{x}_k = - \Grad(\mathbf{x}_k)    
\end{equation}
where \( \Hess(\mathbf{x}_k) \) is the Hessian and \( \mathbf{g}(\mathbf{x}_k) \) is the gradient. Direct inversion of \( \Hess(\mathbf{x}_k) \) is avoided for numerical stability. If \(\Hess(\mathbf{x}_k) \) is singular or ill-conditioned, alternative approaches such as Levenberg-Marquardt regularization, modified Cholesky decomposition, or quasi-Newton methods (e.g., BFGS) can be employed \cite{smejkal_phase_2017, nichita_robustness_2023}. However, no such issues were encountered in our test cases. For implementation, we use Newton's method from \texttt{NLsolve.jl} in Julia, with third-order backtracking Line Search and Trust Region. The gradients and hessian can be computed using automatic differentiation (AD). However, we provide the derivations of the gradients for our approach in \ref{sec:gradAndHessNew} along with the outline of the Hessian matrix, as we intend to use these gradients (of entropy and the constraint function) to compute the Lagrange multiplier, which is further discussed in the following section.

\subsection{Computation of the Lagrange Multiplier}  
In this section, we discuss the computation of the Lagrange multiplier \( \lambda \).
Expanding the stationarity condition \eqref{eqn:kkt1}, we get

\begin{subequations}
    \begin{align}
        \pd{S_{\text{red}}}{N^{(k)}_1} &= -\lambda \pd{\cons}{N^{(k)}_1}, \quad \dots, \quad  \pd{S_{\text{red}}}{N^{(k)}_n} = -\lambda \pd{\cons}{N^{(k)}_n} , \\
        \pd{S_{\text{red}}}{V^{(k)}} &= -\lambda \pd{\cons}{V^{(k)}} , \\
        \pd{S_{\text{red}}}{T} &= -\lambda \pd{\cons}{T} \label{eqn:dSdT_dCdT}.
    \end{align}
\end{subequations}
From equation \eqref{eqn:dSdT_dCdT}, we isolate \(\lambda\) as:

\begin{equation} \label{eqn:lambda_intermediate}
    \lambda = -\frac{\partial S_{\text{red}} / \partial T}{\partial \cons / \partial T}    
\end{equation}
Substituting the expressions from equations \eqref{eqn:dsdT} and \eqref{eqn:dCdT} into equation \eqref{eqn:lambda_intermediate}, we get 
\begin{equation} \label{eqn:lambda}
    \lambda = -\frac{1}{T}.
\end{equation}

The explicit dependence of the Lagrange multiplier \(\lambda\) on temperature \(T\) removes the need to treat \(\lambda\) as an independent optimization variable. This simplification reduces the dimensionality of the problem, as \(\lambda\) is no longer an unknown but is instead directly determined by \(T\). By substituting \eqref{eqn:lambda} into the stationarity condition \eqref{eqn:kkt1}, the optimization process becomes more efficient, as we discuss in detail in the following section.  

\subsection{Lagrangian Simplification} \label{subsec:lagrangian_simplified}
With this choice of the Lagrange multiplier \(\lambda\), the Lagrangian is defined as
\begin{align} \label{eqn:lagrangian_red}
\Lag(\x ) = S_{\text{red}}(\x) - \frac{1}{T} \cons(\x),
\end{align}
where \(\x\) is given by \eqref{eqn:optimisers}. Substituting the expressions for the reduced entropy \(S_{\text{red}}(\x)\) from \eqref{eqn:entropy_red_novel} and the constraint \(\cons(\x)\) from \eqref{eqn:cons_red_novel}, we get
\begin{align}  \label{eqn:lagrangian_S_with_C}
\boxed{ 
\begin{aligned}
\Lag(\x ) &= \sum_{k=1}^{p-1} S\left(T, V^{(k)}, \NBold{k}\right) + S\left(T, V^\star - \sum_{k=1}^{p-1} V^k, \NDiff \right) \\ 
 &\quad - \frac{1}{T} \left( \sum_{k=1}^{p-1} U\left(T, V^{(k)}, \NBold{k}\right) + U\left(T, V^\star - \sum_{k=1}^{p-1} V^k, \NDiff \right) - U^\star \right).
 \end{aligned}
}
\end{align}
We refer to this form of the Lagrangian as the \textbf{Entropy based  Lagrangian} (\(\SCL\) for short) throughout the paper.  Rearranging terms and combining the entropy and internal energy contributions, we obtain
\begin{align} \label{eqn:lagrangian_step2}
\Lag(\x ) &= \sum_{k=1}^{p} \left[ S\left(T, V^{(k)}, \NBold{k}\right) - \frac{U\left(T, V^{(k)}, \NBold{k}\right)}{T}\right] + \frac{U^\star}{T}.
\end{align}
Next, recalling the thermodynamic relation \(A = U - TS\), where \(A\) is the Helmholtz free energy, we get
\begin{align} \label{eqn:lagrangian_step3}
\Lag(\x ) &= \sum_{k=1}^{p} \left[ S\left(T, V^{(k)}, \NBold{k}\right) - \frac{A\left(T, V^{(k)}, \NBold{k}\right) + T S\left(T, V^{(k)}, \NBold{k}\right)}{T}\right] + \frac{U^\star}{T}.
\end{align}
Upon simplifying the terms involving entropy and Helmholtz energy, we arrive at
\begin{align} \label{eqn:lagrangian_step4}
\Lag(\x ) &= -\sum_{k=1}^{p} \frac{A\left(T, V^{(k)}, \NBold{k}\right)}{T} + \frac{U^\star}{T}.
\end{align}
This can be further simplified to
\begin{align} \label{eqn:lagrangian_step5}
\Lag(\x ) &= \frac{U^\star}{T} - \sum_{k=1}^{p} \frac{A\left(T, V^{(k)}, \NBold{k}\right)}{T}.
\end{align}
Finally, segregating the residual terms corresponding to the \(p^{\text{th}}\) phase, we get
\begin{align} \label{eqn:lagrangian_helmholtz}
\boxed{
\begin{aligned}
\Lag(\x ) &= \frac{U^\star - \left( \sum_{k=1}^{p-1} A\left(T, V^{(k)}, \NBold{k}\right) + A\left(T, V^\star - \sum_{k=1}^{p-1} V^k, \NDiff \right)\right)}{T}.
\end{aligned}
}
\end{align}
We refer to this form of the Lagrangian as the \textbf{Helmholtz energy-based Lagrangian} (\(\ACL\) for short). Both equations \eqref{eqn:lagrangian_S_with_C} and \eqref{eqn:lagrangian_helmholtz} represent our novel contribution and provide equivalent formulations for flash calculations. In \eqref{eqn:uncons_entropy_grad}, the Lagrange multiplier \(\lambda\) is treated as an additional unknown, whereas in \(\SCL\) and \(\ACL\) approaches, \(\lambda\) is explicitly determined, reducing the number of unknowns by one. Additionally, \(\SCL\) approach requires more function evaluations compared to \(\ACL\). Specifically, \(\SCL\) involves evaluating both the entropy and the internal energy of each phase, while \(\ACL\) requires only the evaluation of the Helmholtz energy. This suggests that the \(\ACL\) formulation is computationally more efficient than \(\SCL\). We would examine this numerically in \myref{sec:results}. To maximize the entropy, the saddle point of the Lagrangian \(\Lag(\x)\) must be found by solving the system of equations \(\nabla \Lag(\x) = 0\). Additionally, the Hessian matrix in the new formulation simplifies significantly as below:
\begin{align} \label{eqn:hess_simplified}
    \Hess(\x) = \left[ \pdd{\Lag(\x)}{\x_i}{\x_j} \right].
\end{align}
This system can be solved using Newton's method, as discussed in \Cref{sec:newton_method}. The thermodynamic consistency of the proposed formulation will be verified in the next subsection.

\subsection{Thermodynamic consistency check of the new formulation} \label{subsec:gradient_analysis}
In this subsection, we verify the consistency of our new formulation in two key aspects: (1) constraint satisfaction and (2) adherence to thermodynamic equilibrium. To demonstrate consistency, we first show that the formulation with Lagrangian defined as per \Cref{eqn:lagrangian_red} inherently satisfies the constraint of the total internal energy.
The gradient of the Lagrangian with respect to temperature \(T\) is given by:
\begin{align}
    \pd{\Lag}{T} &= \pd{S_{\text{red}}}{T} - \frac{1}{T} \pd{\cons}{T} + \frac{\cons}{T^2} \nonumber \\
   \end{align}
Substituting expressions from \eqref{eqn:dsdT} and \eqref{eqn:dCdT}, we get
\begin{align} \label{eqn:dLdT}
                     \pd{\Lag}{T} &= \frac{1}{T} \sum_{k=1}^{p} C_v\left(T, V^{(k)}, \mathbf{N}^{(k)}\right) - \frac{1}{T} \sum_{k=1}^{p} C_v\left(T, V^{(k)}, 
                     \mathbf{N}^{(k)}\right) + \frac{\cons}{T^2} \nonumber \\
                      &= \frac{\cons}{T^2},
\end{align}
where \(C_v\) represents the heat capacity at constant volume. Using the optimality condition, \(\pd{\Lag}{T} = 0\) yields:
\begin{align}
    \frac{\cons}{T^2} = 0 \implies \cons = 0.
\end{align}
This condition ensures that the constraint \(\cons = 0\) is automatically satisfied, thereby validating the consistency of the formulation.
Next, we verify consistency with respect to thermodynamic equilibrium by computing the gradinets wit respect to the volume and the mole numbers. The gradient of the Lagrangian with respect to the volume \(V^{(k)}\):

\begin{align}
    \pd{\Lag}{V^{(k)}} = \pd{S_{\text{red}}}{V^{(k)}} - \frac{1}{T} \pd{\cons}{V^{(k)}}. 
\end{align}
Substituting the expressions from equations \eqref{eqn:nabla_S_red} and \eqref{eqn:cons_grad}, we obtain
\begin{align*}
    \pd{\Lag}{V^{(k)}} = \frac{\partial P^{(k)}}{\partial T} - \frac{\partial P^{\Res}}{\partial T} 
    - \frac{1}{T}\left( T \left( \pd{P^{(k)}}{T} \right)_{V^{(k)}, \NBoldNoSuper} - P^{(k)} 
    - \left(T \left( \pd{P^{\Res}}{T} \right)_{V^{\Res}, \NBoldNoSuper} - P^{\Res} \right)\right). 
\end{align*}
Here the superscript \(\Res\) denotes evaluation at \(\x^\Res\). After simplification, this reduces to
\begin{align} \label{eqn:dLdV_k}
    \pd{\Lag}{V^{(k)}} = \frac{P^{(k)}}{T} - \frac{P^{\Res}}{T}.
\end{align}
Similarly, the gradient of the Lagrangian with respect to the mole number \(N_1^{(k)}\) is:
\begin{align}
\pd{\Lag}{N_1^{(k)}} &= -\pd{\mu^{(k)}_1}{T} + \pd{\mu_1^{\Res}}{T} 
 - \frac{1}{T} \left(\mu^{(k)}_1 - T \pd{\mu^{(k)}_1}{T} - \left( \mu_1^{\Res} - T \pd{\mu_1^{\Res}}{T} \right)\right) \nonumber \\
 &= -\frac{\mu^{(k)}_1}{T} + \frac{\mu_1^{\Res}}{T}.
\end{align}
Combining these results, the full gradient of the Lagrangian \(\nabla_{\x} \Lag(\x)\) is:
\begin{align} \label{eqn:cons_grad}    
\nabla_{\x} \Lag({\x}) = \begin{pmatrix}
\renderStyle \nabla_{\x} \Lag^{(1)} \\\\ 
\vdots \\\\ 
\renderStyle \nabla_{\x} \Lag^{(p-2)} \\\\
\renderStyle \nabla_{\x} \Lag^{(p-1)} \\\\
\renderStyle \pd{\Lag}{T}
\end{pmatrix},
\end{align}
where the individual entry 
\begin{align}
    \nabla_{\x}{\Lag}^{(k)} = \begin{pmatrix}  
          \renderStyle -\frac{\mu^{(k)}_1}{T} + \frac{\mu_1^{\Res}}{T} \\\\ 
           \vdots \\\\
          \renderStyle -\frac{\mu^{(k)}_n}{T} + \frac{\mu^{\Res}_n}{T} \\\\ 
          \renderStyle \frac{P^{(k)}}{T} - \frac{P^{\Res}}{T}
    \end{pmatrix}.
\end{align}
 and \(\pd{\Lag}{T}\) is given by \myeqref{eqn:dLdT}.
The final gradients of the Lagrangian are structurally identical to those reported by Smejkal et al. \cite{smejkal_phase_2017}, with the key distinction that our formulation allows all functions to be evaluated directly as a function of \(T, V\) and \(\NBoldNoSuper\), whereas, Smejkal's formulation requires an inner Newton iteration to first determine the temperature. The optimality condition \(\nabla_{\x} \Lag({\x}) = 0\) leads to the following system of equations.
\begin{subequations}
    \begin{align}
        \mu^{(1)}_1 = \mu^{(2)}_1 &= \dots = \mu^{\Res}_1 \\
        \mu^{(1)}_2 = \mu^{(2)}_2 &= \dots = \mu^{\Res}_2 \\ 
        &\vdots \\
        \mu^{(1)}_n = \mu^{(2)}_n &= \dots = \mu^{\Res}_n \\
        P^{(1)} = P^{(2)} &= \dots = P^{\Res}
    \end{align}
\end{subequations}
i.e. the chemical potential of each component is equal across all coexisting phases and the pressure of each phase is equal. This is consistent with the principles of thermodynamic equilibrium. 

With the theoretical framework established, we now proceed to the results section, where we present numerical results obtained using the proposed methodology which leverages the Lagrangians defined in Equations \eqref{eqn:lagrangian_S_with_C} and \eqref{eqn:lagrangian_helmholtz}. 



\section{Results} \label{sec:results}
In this section, we present the results obtained using our proposed approach and compare them with existing literature. Our treatment focuses exclusively on the two-phase test cases examined by Castier~\cite{castier_solution_2009}, Smejkal et al.~\cite{smejkal_phase_2017}, and Bi et al.~\cite{bi_efficient_2020}. These problems have also been discussed by Nichita~\cite{nichita_robustness_2023} in the context of VT stability analysis. Specifically, we consider Problems 1--6 from these studies, along with a pure component test case introduced by Smejkal. These problems are defined in \myref{tab:problem_1_4_spec}, \ref{tab:problem_5_6_spec} and \ref{tab:problem_pure_co2_spec}. 
Notably, no variable scaling was employed during the optimization process. In contrast, Smejkal et al. did not explicitly state whether variable scaling was used in their approach.

We begin by discussing the outcomes of the stability analysis, which serve as the foundation for determining the initial phase split. These results are then used to perform flash calculations, the details of which are presented subsequently. First, we validate our results with literature followed by a discussion of the speedup gains. Finally, we conclude this section with a comparison of the two forms of Lagrangian derived in this work.

\begin{table} [H]
    \centering
    \begin{tabular}{lrrrr}
    \hline
        Property & Problem 1  & Problem 2 & Problem 3 & Problem 4\\
        \hline
        $U$[J] & -756500.8 & -1511407.6 & -331083.7 & -636468 \\
        $V$[\volUnit] & 52869 & 4268.1 & 80258.1 & 9926.71 \\
        $N_{c_1}$[mol] & 10 & 0.95 & 15.1 & 10 \\         
        $N_{H_2S}$[mol] & 90 & 99.05 & 84.9 & 90 \\     
        \hline
    \end{tabular}
    \caption{Specification : Problems 1--4}
    \label{tab:problem_1_4_spec}
\end{table}
\begin{table}[H]
    \centering
    \begin{tabular}{lrr}
        \hline
        Property  & Problem 5 & Problem 6 \\
        \hline
        $U$ [J]           & -16272506.4  & 24858.2 \\
        $V$ [\volUnit]      & 479845       & 289380.3 \\
        $N_{C_2}$ [\concUnit] & 10.8         & 10.8 \\
        $N_{C_3H_6}$ [\concUnit] & 360.8     & 360.8 \\
        $N_{C_3}$ [\concUnit] & 146.5        & 146.5 \\
        $N_{iC_4}$ [\concUnit] & 233         & 233 \\
        $N_{nC_4}$ [\concUnit] & 233         & 233 \\
        $N_{C_5}$ [\concUnit] & 15.9         & 15.9 \\
        \hline
    \end{tabular}
    \caption{Specification : Problems 5-6}
    \label{tab:problem_5_6_spec}
\end{table}

\begin{table}[H]
    \centering
    \begin{tabular}{cccc}
    \hline
        Property & $U$ [J] & $V$ [m$^3$] & $N_{\text{CO}_2}$ [mol] \\ \hline
        Value & -87211375.744478 & 1 & 10000 \\ \hline
    \end{tabular}
    \caption{Specification: Pure component CO$_2$}
    \label{tab:problem_pure_co2_spec}
\end{table}

For all calculations, the Peng-Robinson equation of state (EOS)~\cite{peng_new_1976}, based on Helmholtz energy, is employed. Additional details regarding this EOS can be found in~\ref{app:PR_EOS}.

\subsection{Stability Analysis} \label{sec:results_stability}
While the primary focus of this paper is on phase split calculations, we first present the results of the stability analysis, as these results provide the initial guesses for the phase split calculations. We have obtained these results using the methodology discussed in \Cref{sec:stability_analysis}. Our study reports the local minimum for each problem, with the results summarized in Tables \ref{tab:problem_1_4_stability_results}, \ref{tab:problem_5_6_stability_results} and \ref{tab:pure_co2_stability_results}.  For each case, we report the computed values of temperature, component concentrations, and the tangent plane distance function \(D\), as defined in \myeqref{eqn:TPD_func}, with the results reported to two significant digits. However, for values smaller than 1, the results are reported to four significant digits. In all cases, our local minima are in close agreement with the values (either global or local) reported by Nichita~\cite{nichita_robustness_2023} for the multicomponent case and Smejkal et al.~\cite{smejkal_phase_2017} for the single component case. The stability analysis reveals minimal discrepancies in concentration values, with errors remaining below 0.085\%. The largest errors occur in Problem 6, with the highest being 0.085\% for \(c^\prime_{C_5}\). In the following section, we utilize these stability results to initialize the phase split calculations. 

\begin{table}[H]
    \centering
    \begin{tabular}{ccccccccc}
        \hline
         & \multicolumn{2}{c}{Problem 1} & \multicolumn{2}{c}{Problem 2} & \multicolumn{2}{c}{Problem 3} & \multicolumn{2}{c}{Problem 4} \\
        \cline{2-9}
        Property & \LitStabRes & \OurStabRes  &  \LitStabRes & \OurStabRes & \LitStabRes & \OurStabRes & \LitStabRes & \OurStabRes \\
        \hline
        $T$ [K]                    &  151.83   & 151.83     & 291.91      & 291.91      & 297.84    & 297.84     & 361.80     & 361.80      \doubleEnter
        $c_{c_1}'$ [\concUnit]     &  104.13  & 104.12     & 146.11      & 146.18      & 188.14    & 188.14     & 1011.37    & 1011.36     \doubleEnter
        $c_{H_2S}'$ [\concUnit]    &  564.39  & 564.35     & 736.15      & 736.58      & 1057.84   & 1057.84    & 10056.7    & 10037.91      \doubleEnter
        $D$ [Pa/K]                 &  875.34   & 875.45     & 26771.1     & 26722       & 0.0       & 2.08e-12   & 0.5063     &  0.467      \\
        \hline
    \end{tabular}
    \caption{Results of stability analysis: Nichita~\cite{nichita_robustness_2023} vs Our Results }
    \label{tab:problem_1_4_stability_results}
\end{table}

\begin{table}[H]
    \centering
    \begin{tabular}{lrrrr}
        \hline
        & \multicolumn{2}{c} {Problem 5} &   \multicolumn{2}{c}{Problem 6} \\
        \cline{2-5}
         Property      &           \LitStabRes   &    \OurStabRes &  \LitStabRes          &    \OurStabRes   \\
        \hline
        $T$[K]                        &  122.97   &    122.97         &   394.54   & 394.54             \doubleEnter
        $c^\prime_{C_2}$ [mol]      & 0.3294    &    0.3294         &   46.41    & 46.41      \doubleEnter
        $c^\prime_{C_3H_6}$ [mol]   & 3.10      &    3.10           &   1739.38  & 1738.53     \doubleEnter
        $c^\prime_{C_3}$ [mol]      & 0.9066    &    0.9066         &   719.16   & 718.79     \doubleEnter
        $c^\prime_{iC_4}$ [mol]     & 0.3860    &    0.3860         &   1262.45  & 1261.59     \doubleEnter
        $c^\prime_{nC_4}$ [mol]     & 0.2934    &    0.2934         &   1305.65  & 1304.69    \doubleEnter
        $c^\prime_{C_5}$ [mol]      & 0.0038    &    0.0038         &   101.09   & 101.00     \doubleEnter
        $D$[Pa/K]                     & 35298.75  &    35298.74       &   16.3045  & 16.10         \\ 
        \hline
    \end{tabular}
    \caption{Results of stability analysis Nichita~\cite{nichita_robustness_2023} vs Our Results}
    \label{tab:problem_5_6_stability_results}
\end{table}

\begin{table}[H]
    \centering
    \begin{tabular}{ccc}
    \hline
        Property & Smejkal & \OurStabRes \\
        \hline
        $T$ [K]                     &   280.0      &   280.0        \\         
        $c^{\prime}$ [\concUnit]    &   19469.17   &   19487.12    \\ 
        $D$ [Pa/K]                  &   4608.22    &   4608.27      \\
        \hline
    \end{tabular}
    \caption{Results for pure CO$_2$ from stability analysis. Smejkal et al.~\cite{smejkal_phase_2017} vs Our Results}
    \label{tab:pure_co2_stability_results}
\end{table}

\subsection{Flash calculations} \label{sec:results_flash}
In this section, we present the initial guesses derived from stability analysis, generated using the  \myref{alg:init_guess} described in \myref{sec:initial_guess_from_stability}. While Smejkal et al. \cite{smejkal_phase_2017} highlight the use of stability analysis to obtain initial guesses for flash calculations, their work does not explicitly provide these values for all the test cases, limiting the reproducibility of their results. To bridge this gap, we report the detailed initial guesses obtained from our stability analysis, followed by the results of the corresponding flash calculations. The initial guesses are comprehensively summarized in Tables \ref{tab:initial_guess_p1234}, \ref{tab:initial_guess_p56} and \ref{tab:co2_initial_guess}, with the results reported to four significant digits.

Using these initial guesses, we perform flash calculations using our proposed approach. For all results presented here, we have used the entropy-based Lagrangian defined as per \Cref{eqn:lagrangian_S_with_C}. The results with the Helmholtz energy-based Lagrangian formulation are essentially the same and are therefore not included here. Tables \ref{tab:flash_1_2}, \ref{tab:flash_3_4} and \ref{tab:flash_5_6} present the results using Netwon's method with a third order backtracking linesearch. The stopping criterion is set to a relative tolerance of \(1 \times 10^{-8}\). The results are reported to six significant digits. In addition to the internal energy, volume and mole numbers, we also report the entropy of the reference phase and the two phase system, denoted as \(S^{I}\) and \(S^{II}\), respectively. A reasonable agreement is observed with the results reported by Smejkal et al. \cite{smejkal_phase_2017} for problems 1-6.

\begin{table}[H]
    \centering
    \resizebox{\textwidth}{!}{%
    \begin{tabular}{lrrrrrrrr}
        \hline
        Property   & \multicolumn{2}{c}{Problem 1} & \multicolumn{2}{c}{Problem 2} & \multicolumn{2}{c}{Problem 3} & \multicolumn{2}{c}{Problem 4} \\
        \cline{2-9}
                 & Phase 1 & Phase 2 & Phase 1 & Phase 2 & Phase 1 & Phase 2 & Phase 1 & Phase 2 \\
        \hline
        $N_{c_1}$[mol]  & 0.0003   & 9.9997   & 0.0195  & 0.9305   & 0.0005  & 15.0995  & 1.2549  & 8.7451  \\
        $N_{H_2S}$[mol] & 28.5304  & 61.4696  & 0.0982  & 98.9518  & 0.0583  & 84.8417  & 12.4787  & 77.5213  \\
        $V$[\volUnit]        & 0.0008  & 0.0520  & 0.0001  & 0.0041  & 2.45e-6  & 0.08026  & 0.0012  & 0.0087  \\
        $U$[J]        & -717694  & -38806   & -379.56  & -1.511e6 & -891.17  & -330193  & -94307.8  & -542160  \\
        \hline
    \end{tabular}
    }
    \caption{Initial guesses obtained from stability analysis for Problems 1, 2, 3, and 4}
    \label{tab:initial_guess_p1234}
\end{table}

\begin{table}[H]
    \centering    
    \begin{tabular}{lrrrr}
        \hline
        Property  & \multicolumn{2}{c}{Problem 5} & \multicolumn{2}{c}{Problem 6} \\
        \cline{2-5}
                         & Phase 1       & Phase 2       & Phase 1       & Phase 2       \\
        \hline
        $N_{C_2}$ [mol]  & 0.6328        & 10.2672       & 0.8394        & 9.9624        \\
        $N_{C_3H_6}$ [mol] & 84.6926      & 276.1074      & 31.4435       & 329.3565      \\
        $N_{C_3}$ [mol]  & 22.4028       & 124.0972      & 13.0002       & 133.4998      \\
        $N_{iC_4}$ [mol] & 26.3615       & 206.6385      & 22.8175       & 210.1825      \\
        $N_{nC_4}$ [mol] & 82.5333       & 150.4667      & 23.5972       & 209.4028      \\
        $N_{C_5}$ [mol]  & 9.1875        & 6.7125        & 1.8268        & 14.0732       \\
        $V$ [\volUnit]   & 0.0150        & 0.4648        & 0.0181        & 0.2713        \\
        $U$ [J]          & -8.2275e6     & -8.0450e6     & -211881.92    & 236740.12     \\
        \hline
    \end{tabular}
    \caption{Initial guesses from stability analysis for Problems 5 and 6}
    \label{tab:initial_guess_p56}
\end{table}

\begin{table}[H]
    \centering
    \begin{tabular}{lc}
        \hline
        Property   & Phase 1\\
        \hline
        $N_{CO_2}$ [mol]    & 2435.89 \\
        $V$ [\volUnit]     & 0.125    \\
        $U$ [J]            & $-3.129 \times 10^7$ \\
        \hline
    \end{tabular}    
    \caption{Initial guesses obtained from stability analysis for pure component (\coo)}
    \label{tab:co2_initial_guess}
\end{table}

\begin{table}[H]
    \centering
    \begin{tabular}{lrrrrr}
        \hline
        & \multicolumn{2}{c}{Problem 1} & \multicolumn{2}{c}{Problem 2} \\
        \cline{2-5}
        & Smejkal & \OurStabRes & Smejkal & \OurStabRes \\
        \hline
        $U$ [J]             & –211544.585681   & -211544.596326   & –1510985.753624   &  -1510985.755666 \\
        $V$ [cm$^3$]        & 51366.638771     & 51366.638597     & 4165.673900       &  4165.674425 \\
        $N_{c_1}$ [mol]     & 9.664320         & 9.664319         & 0.930730          &  0.930730 \\
        $N_{H_2S}$ [mol]    & 54.315978        & 54.315976        & 98.941685         &  98.941685 \\
        \hline
        $S^{I}$ [J K$^{-1}$] & –4847.824318 & -4847.824867 & –7391.709463 & -7391.709647 \\
        $S^{II}$ [J K$^{-1}$] & –4335.499136 & -4335.499558 & –7390.326639 & -7390.326837 \\
        \hline
    \end{tabular}
    \caption{Comparison of Flash Results for Problems 1 and 2: Smejkal et al.~\cite{smejkal_phase_2017} vs. Our Results}
    \label{tab:flash_1_2}
\end{table}

\begin{table}[H]
    \centering
    \begin{tabular}{lrrrrr}
        \hline
        & \multicolumn{2}{c}{Problem 3} & \multicolumn{2}{c}{Problem 4} \\
        \cline{2-5}
        & Smejkal & \OurStabRes & Smejkal & \OurStabRes \\
        \hline
        $U$ [J]             & –330516.922985   &  -330516.953672    & –390660.034825 & -390689.64236 \doubleEnter
        $V$ [cm$^3$]        & 80256.537494     &  80256.537579      & 6414.083981    & 6414.415486 \doubleEnter
        $N_{c_1}$ [mol]     & 15.099651        &  15.099651         & 6.448582       & 6.448928 \doubleEnter
        $N_{H_2S}$ [mol]    & 84.862887        &  84.862889         & 56.390527      & 56.394270 \doubleEnter
        \hline
        $S^I$ [J K$^{-1}$] & –2613.988230 & -2613.988418 & –4579.402758 & -4579.403289 \doubleEnter
        $S^{II}$ [J K$^{-1}$] & –2613.987835 & -2613.988023 & –4579.402147 & -4579.402679\doubleEnter
        \hline
    \end{tabular}
    \caption{Comparison of Flash Results for Problems 3 and 4: Smejkal et al.~\cite{smejkal_phase_2017} vs. Our Results}
    \label{tab:flash_3_4}
\end{table}

\begin{table}[H]
    \centering
    \begin{tabular}{lrrrrr}
        \hline
        & \multicolumn{2}{c}{Problem 5} & \multicolumn{2}{c}{Problem 6} \\
        \cline{2-5}
        & Smejkal & \OurStabRes & Smejkal & \OurStabRes \\
        \hline
        $U$ [J]             & –379886.931385       &  -380012.963119    & 174870.975415      &  174842.436972  \\
        $V$ [cm$^3$]        & 401197.390420        &  401192.630291     & 273147.423428      &  273150.189814   \\
        $N_{C_2}$ [mol]     & 4.203436             &  4.242459          & 10.064693          &  10.066498       \\
        $N_{C_3H_6}$ [mol]  & 68.225832            &  68.231202         & 333.710698         &  333.715455        \\
        $N_{C_3}$ [mol]     & 24.416960            &  24.419097         & 135.325654         &  135.327702     \\
        $N_{iC_4}$ [mol]    & 18.529159            &  18.531724         & 213.665513         &  213.668936       \\
        $N_{nC_4}$ [mol]    & 13.885437            &  13.887650         & 213.118914         &  213.122442         \\
        $N_{C_5}$ [mol]     & 0.325600             &  0.325674          & 14.391190          &  14.391459       \\
        \hline
        $S^{I}$ [J K$^{-1}$] & –73647.697512 &  -73640.643944 & –9052.552759 & -9052.541673 \\
        $S^{II}$ [J K$^{-1}$] & –54939.068244 & -54937.804163 & –9052.431373 & -9052.420341 \\
        \hline
    \end{tabular}
    \caption{Comparison of Flash Results for Problems 5 and 6: Smejkal et al.~\cite{smejkal_phase_2017} vs. Our Results}
    \label{tab:flash_5_6}
\end{table}

 To further evaluate the generality and robustness of our method, we also consider a single-component test case, as discussed by Smejkal et al. \cite{smejkal_phase_2017}, with specifications defined in  \myref{tab:problem_pure_co2_spec}. The stability analysis (see \myref{tab:pure_co2_stability_results}) reveals that the fluid is unstable as a single-phase fluid. Based on this analysis, an initial phase split was obtained, as shown in Table \ref{tab:co2_initial_guess}.  Flash calculations are subsequently performed using this initial phase split, and the results are presented in Table \ref{tab:flash_co2_single_phase}. Our findings show excellent agreement with the results reported in the literature \cite{smejkal_phase_2017}. Following this validation, we proceed to evaluate the computational efficiency of our method, comparing the speedup relative to the approach of Smejkal et al.

\begin{table}[H]
    \centering
    \begin{tabular}{lrr}
        \hline
        Property   &Smejkal & \OurStabRes \\
        \hline
        $U$ [J]         & -16873789.390417    & -16873791.656255 \\
        $V$ [\volUnit]       & 481283.619636  & 481283.486064 \\        
        $N_{CO_2}$ [mol]     & 2818.038884 & 2818.038719 \\
        \hline
        $S^I$ [J/K]       &  -584388.217059 & -584388.23982 \\
        $S^{II}$ [J/K]     & -583476.321606 & -583476.346351 \\
        \hline
    \end{tabular}
    \caption{Comparison of Flash Results for pure CO$_2$: Smejkal vs. Our Results}
    \label{tab:flash_co2_single_phase}
\end{table}

\subsubsection{Speedup and Robustness} \label{sec:speedup}

We now turn our attention to the computational speedup achieved by our TVN approach compared to the UVN formulation. Both formulations are compared by directly using the same nonlinear solver in \texttt{Julia} without employing any variable scaling.  The results,  obtained with a relative tolerance of \(1 \times 10^{-6}\), are summarized in Tables \ref{tab:iterations_combined} and \ref{tab:execution_time_combined_scaled}.

The TVN approach consistently requires fewer or equal outer iterations compared to the UVN method across all test problems. For instance, in Problem P1, both methods require 9 outer iterations. In contrast, in Problem P4, the TVN approach converges with 7 outer iterations, whereas the UVN approach fails to converge when using Line Search. However, when Trust Region method is used, the UVN formulation requires 15 outer iterations, whereas the TVN formulation requires only 7 iterations for the same problem. This demonstrates the robustness of the TVN formulation, which reliably converges for all problems using both Line Search and Trust Region methods. 

A key advantage of the TVN approach is the elimination of inner iterations, which are computational bottleneck in the UVN method. While UVN requires 20 to 987 inner iterations, TVN bypasses this overhead entirely, significantly reducing computational costs in the current implementation. This difference is evident in the single-component test case (P\coo), where TVN completes with 32 outer iterations (Line Search) and 117 (Trust Region), compared to UVN’s 77 and 175 outer iterations plus 964–987 inner iterations. Clearly this is reflected in the execution time, for example for Problem P1, TVN requires 0.52 ms (Line Search) and 0.41 ms (Trust Region), whereas UVN takes 3.15 ms and 2.78 ms, respectively. In the P\coo case, TVN achieves 0.69 ms (Line Search) and 1.18 ms (Trust Region), compared to UVN’s 263.24 ms and 481.01 ms. This performance underscores the scalability of the TVN approach for practical applications.

\begin{table}[H]
    \centering
    \begin{tabular}{l|ccc|ccc}
        \cline{2-7}
        & \multicolumn{3}{c|}{\textbf{Line Search}} & \multicolumn{3}{c}{\textbf{Trust Region}} \\
        \hline
       & \textbf{TVN} & \multicolumn{2}{c|}{\textbf{UVN}} & \textbf{TVN} & \multicolumn{2}{c}{\textbf{UVN}} \\
        \cline{3-4} \cline{6-7}
        Problem & & Outer & Inner & & Outer & Inner \\
        \hline
        P1   & 9  & 9   & 65  & 9   & 9   & 64  \\
        P2   & 4  & 4   & 20  & 4   & 4   & 20  \\
        P3   & 4  & 4   & 30  & 4   & 5   & 30  \\
        P4   & 7  & x   & x   & 7   & 15  & 81  \\
        P5   & 10 & 10  & 71  & 10  & 10  & 68  \\
        P6   & 4  & 5   & 60  & 5   & 5   & 60  \\
        P\coo & 32  & 77 & 964 & 117 & 175 & 987 \\
        \hline
    \end{tabular}
    \caption{Iteration counts for TVN ($\SCL$) and UVN formulations. Cases that did not converge are marked with x}
    \label{tab:iterations_combined}
\end{table}

\begin{table}[H]
\centering
\begin{tabular}{l|ccc|ccc}
\cline{2-7}
        & \multicolumn{2}{c}{\textbf{Line Search}} & &\multicolumn{2}{c}{\textbf{Trust Region}} \\
        \cline{1-7}
        & &\multicolumn{2}{c|}{\textbf{UVN}} & &\multicolumn{2}{c}{\textbf{UVN}} \\
        \cline{3-4} \cline{6-7} 
        Problem & TVN (ms) & No Scale (ms) & Scaled (ms) & TVN (ms) & No Scale (ms) &  Scaled(ms) \\
        \hline
P1   & 0.52  & 3.15   & 3.66  & 0.41   & 2.78   & 3.67   \\
P2   & 0.28  & 1.86   & 2.09  & 0.22   & 1.22   & 1.33  \\
P3   & 0.17  & 1.09   & 1.23  & 0.22   & 1.48   & 2.08  \\
P4   & 0.47  & x      & 12.14 & 0.38   & 3.66   & 4.23  \\
P5   & 4.59  & 42.47  & 40.3  & 5.29   & 42.86  & 41.48 \\
P6   & 2.37  & 26.67  & 26.73 & 3.16   & 32.54  & 31.46 \\
P\coo & 0.69 & 263.24 & 45.69 & 1.18  & 481.01 &  1.86  \\
\hline
\end{tabular}
\caption{Execution time comparison (in milliseconds) for TVN ($\SCL$) and UVN formulations. The execution time reported is the average calculation time after solving each problem 100 times. Cases that did not converge are marked with x.}
\label{tab:execution_time_combined_scaled}
\end{table}

It is worthwhile noting that the convergence challenges observed in the UVN method under Line Search can be mitigated by applying variable scaling, where all variables are normalized by their respective total specified quantities. This approach has been tested and confirmed effective, as illustrated in \myref{tab:execution_time_combined_scaled}. For instance, in Problem P\coo, scaling reduces UVN’s execution time from 263.24 ms to 45.69 ms with Line Search and from 481.01 ms to 1.86 ms with the Trust Region method. These results emphasize the importance of variable scaling in improving the performance and robustness of the UVN method.

We further point out that the TVN approach benefits from having better-scaled variables from the start, for example, temperature typically varies within a narrow range of a few hundred kelvins, unlike internal energy, which spans a wide range of large negative and positive values. This enhances the numerical stability and contributes to the robustness of the TVN approach. Overall, the robustness combined with reduced computational time positions of the TVN approach establishes it as a superior alternative to the conventional UVN method for thermodynamic flash calculations.

\subsubsection{Comparison of Two Lagrangian Formulations}
This section compares two novel Lagrangian formulations, defined in equations \eqref{eqn:lagrangian_S_with_C} and \eqref{eqn:lagrangian_helmholtz}, using the iteration counts and the execution times presented in Tables~\ref{tab:tvn_lagrangian_comparison_iterations} and  \ref{tab:tvn_lagrangian_comparison_times} respectively. The first formulation, which employs entropy based Lagrangian ($\SCL$), while the second formulation, which utilizes the Helmholtz energy-based Lagrangian ($\ACL$).

The iteration counts for both approaches are similar, however, the $\ACL$ formulation consistently outperforms the $\SCL$ formulation in execution time across all problems and methods (Line Search and Trust Region). This efficiency is attributed to the fewer function evaluations required by \(\ACL\) approach, as discussed in \ref{subsec:lagrangian_simplified}. For instance, in P5 using Newton with Line Search, $\ACL$ takes only 1.32 ms whereas \(\SCL\) takes 4.59 ms. Additionally, the Trust Region method generally exhibits longer execution times compared to the Line Search method, particularly for the \(\SCL\) formulation, as observed in problems P5 and P6. It is important to note that the execution time comparisons for problems P1\textendash P4 are less meaningful, as these problems already require minimal computational time, rendering execution time an unreliable metric these specific cases. Overall, these results demonstrate the computational efficiency of the Helmholtz-based Lagrangian and highlight the significant impact of optimization method selection (Line Search vs. Trust Region) on performance. 

\begin{table}[H]
\centering
\begin{tabular}{lcc|cc}
\hline
& \multicolumn{2}{c}{\textbf{Line Search}} & \multicolumn{2}{|c}{\textbf{Trust Region}} \\
\cline{2-5}
Problem & $\SCL$ (ms) & $\ACL$ (ms) & $\SCL$ (ms) & $\ACL$ (ms) \\
\hline
P1      & 9     & 10  & 9     & 10    \\
P2      & 4     & 4   & 4     & 4     \\
P3      & 4     & 4   & 4     & 4     \\
P4      & 7     & 8   & 7     & 8     \\
P5      & 10    & 10  & 10    & 10    \\
P6      & 4     & 5   & 5     & 5     \\
P\coo   & 32    & 33  & 117   & 117   \\
\hline
\end{tabular}
\caption{Iteration counts comparison for $\SCL$ and $\ACL$ formulations.}
\label{tab:tvn_lagrangian_comparison_iterations}
\end{table}

\begin{table}[H]
\centering
\begin{tabular}{lcc|cc}
\hline
& \multicolumn{2}{c}{\textbf{Line Search}} & \multicolumn{2}{|c}{\textbf{Trust Region}} \\
\cline{2-5}
Problem & $\SCL$ (ms) & $\ACL$ (ms) & $\SCL$ (ms) & $\ACL$ (ms) \\
\hline
P1   & 0.52  & 0.24  & 0.41  & 0.21 \\
P2   & 0.28  & 0.15  & 0.22  & 0.13 \\
P3   & 0.17  & 0.11  & 0.22  & 0.14 \\
P4   & 0.47  & 0.23  & 0.38  & 0.23 \\
P5   & 4.59  & 1.32  & 5.29  & 2.53 \\
P6   & 2.37  & 0.7   & 3.16  & 0.68 \\
P\coo & 0.69  & 0.51  & 1.18  & 0.92 \\
\hline
\end{tabular}
\caption{Execution time comparison (in milliseconds) for $\SCL$ and $\ACL$ formulations. The reported calculation time is the average time after repeating 100 calculations for each problem.}
\label{tab:tvn_lagrangian_comparison_times}
\end{table}

\section{Conclusion} \label{sec:conclusion}  
In this work, we propose a novel reformulation of the UVN flash problem that eliminates the need to invert the equation of state (EOS) during the optimization process, significantly improving the efficiency of flash calculations. By transitioning from the unconstrained UVN space to the constrained TVN space, we simplify the numerical approach, enhancing its robustness and facilitating faster convergence. This TVN reformulation addresses a key limitation in previous methods, which requires costly EOS inversions as part of the optimization procedure.

We applied the method of Lagrange multipliers to transform the constrained optimization problem into an unconstrained one.  By deriving the necessary gradients and Hessian, we obtained an explicit expression for the Lagrange multiplier in terms of temperature, eliminating the need to treat it as an independent variable. This simplification led to two Lagrangian formulations, namely, the entropy-based Lagrangian (\(\SCL\)) and the Helmholtz energy-based Lagrangian (\(\ACL\)). 

We also provided an explicit algorithm for generating high-quality initial guesses directly from stability analysis results. This critical step greatly facilitates the convergence of the flash algorithm. 
We subsequently applied our novel reformulation of the UVN flash and validated the results against literature data. To evaluate the performance of our method, we employed two variants of Newton's: Line Search and Trust Region. Our TVN reformulation demonstrated successful convergence with both approaches, highlighting its robustness and reliability. We tested our approach on a single-component system, that further demonstrates its wide applicability. 

Finally, a comparison of our entropy-based Lagrangian formulation against the traditional UVN approach was performed. The computational results show significant performance improvements, including: avoiding the need to perform EOS inversions and a reduction in the number of outer loop iterations. These improvements led to a reduction in execution time without sacrificing accuracy. Additionally, both of the proposed Lagrangian formulations were compared in terms of execution time. The results indicate that the Helmholtz-based formulation consistently outperforms the entropy-based formulation in computational efficiency across all test cases.

\section*{CRediT author statement}


\textbf{Pardeep Kumar}:
Conceptualization,
Methodology,
Software,
Validation,
Analysis,
Writing -- Original Draft

\textbf{Patricio I. Rosen Esquivel}:
Project Administration,
Funding Acquisition,
Conceptualization,
Supervision,
Writing - Review \& Editing

\section*{Declaration of Generative AI and AI-assisted technologies in the writing process}


During the preparation of this work the authors used GitHub Copilot in order to
propose wordings and mathematical typesetting. After using this tool/service,
the authors reviewed and edited the content as needed. The authors take full
responsibility for the content of the publication.
\section*{Declaration of competing interest}


The authors declare that they have no known competing financial interests or
personal relationships that could have appeared to influence the work reported
in this paper.
\section*{Acknowledgments}  
This research was generously supported by Shell Projects and Technology, and we deeply appreciate their invaluable contribution.

We would like to express our sincere gratitude to \textbf{Prof. Ruud Henkes} (TU Delft) and \textbf{Prof. Benjamin Sanderse} (CWI Amsterdam, TU Eindhoven) for their expert guidance and insightful contributions throughout this research.

Our heartfelt thanks also go to \textbf{Dr. Jannis Teunissen} (CWI Amsterdam) and \textbf{Dr. Marius Kurz} (CWI Amsterdam) for their constructive feedback on the manuscript.

Finally, we are especially grateful to \textbf{Prof. Jiri Mikyska} (Czech Technical University in Prague) for his invaluable insights and stimulating discussions on UVN Flash.

\bibliographystyle{abbrv}

\appendix

\section{Gradient and Hessian Computation in new Formulation} \label{sec:gradAndHessNew}
In this section, we present the methodology for computing the gradient and Hessian of the entropy and the constraint functions in the proposed formulation. These derivatives are critical for solving the optimization problem using Newton's method, as outlined in the preceding section. The discussion is structured as follows: we first present the computation of the gradient, followed by the derivation of the Hessian.

\subsection{Gradients computation}
In this subsection, we discuss the evaluation of the gradients of the the entropy function defined by \myeqref{eqn:entropy_red_novel} and the constraint function defined by \myeqref{eqn:cons_red_novel}. Using these gradients, we then compute the value of the Lagrange multiplier \(\lambda\). The gradient \(\nabla S_{\text{red}}(\xNew)\) of the reduced entropy function \(S_{\text{red}}(\xNew)\) is defined as:


\begin{align} \label{eqn:entropy_grad_new}    
\nabla S_{\text{red}}({\xNew}) = \begin{pmatrix}
\nabla S^{(1)}_{\text{red}} \\ 
\vdots \\ 
\nabla S^{(p-2)}_{\text{red}} \\\\
\nabla S^{(p-1)}_{\text{red}} \\\\
\renderStyle \pd{S_{\text{red}}}{T}
\end{pmatrix}
\end{align}
$\forall k \in \{1, \dots, p-1\}$, \quad \(\nabla S_{\text{red}}^{(k)} \in \mathbb{R}^{n+1}\),  and \(\pd{S_{\text{red}}}{T} = \sum_{k=1}^{p} \pd{S^{(k)}}{T} \in \mathbb{R}\), where \(S^{(k)} = S\left(T, V^{(k)}, \NBold{k}\right)\) is the entropy of the phase $k$. The individual entries of \(\nabla S_{\text{red}}^{(k)}\) are given as below.

\begin{align}
    \nabla S^{(k)}_{\text{red}} = \colvec{
        \renderStyle \pd{S_{\text{red}}}{N^{(k)}} \\\\
        \vdots \\\\
        \renderStyle \pd{S_{\text{red}}}{N^{(n)}} \\\\
        \renderStyle \pd{S_{\text{red}}}{V^{(k)}} 
        } = \colvec{
       \renderStyle  \pd{S^{(k)}_{\text{red}}}{N^{(k)}_1} - \pd{S^{\Res}_{\text{red}}}{N^{(k)}_1} \\\\
        \vdots \\\\
       \renderStyle  \pd{S^{(k)}_{\text{red}}}{N^{(k)}_n} - \pd{S^{\Res}_{\text{red}}}{N^{(k)}_n} \\\\
       \renderStyle  \pd{S^{(k)}_{\text{red}}}{V^{(k)}} - \pd{S^{\Res}_{\text{red}}}{V^{(k)}} 
        }.
\end{align}
We can simplify the partial derivatives using thermodynamic identities as follows:
\begin{align} 
    \pd{S_{\text{red}}}{T} &= \sum_{k=1}^{p-1} \pd{S\left(T, V^{(k)}, \mathbf{N}^{(k)}\right)}{T} + \pd{S\left(T, V^\star - \sum_{k=1}^{p-1} V^{(k)}, \mathbf{N}^{\Res} \right)}{T} \nonumber \\
    &= \frac{1}{T} \sum_{k=1}^{p} C_v\left(T, V^{(k)}, \mathbf{N}^{(k)}\right). \label{eqn:dsdT}
\end{align}
Furthermore, the thermodynamic identity for the volume derivative of entropy is given by:
\begin{align}
    \left(\pd{S}{V}\right)_{T, \mathbf{N}} = \left(\pd{P}{T}\right)_{V, \mathbf{N}}.
\end{align}
Next, for the derivative with respect to \( N \), we can substitute \(S\) in terms of Helmholtz energy \(A\) as follows:
\begin{align}
    \left(\pd{S}{N}\right)_{T, V} = \left(\pd{\left(-\pd{A}{T}\right)_{V, N}}{N}\right)_{T, V}.
\end{align}
Since \( V \) is constant, we consider only \( T \) and \( N \) as variables, yielding:
\begin{align*}
    \left(\pd{S}{N}\right)_{T} &= \left(\pd{\left(-\pd{A}{T}\right)_{N}}{N}\right)_{T} \\
     & = -\pd{\mu}{T}.
\end{align*}
Finally, the gradient of the reduced entropy for phase \( k \) is given by:
\begin{align} \label{eqn:nabla_S_red}
    \nabla S^{(k)}_{\text{red}} = \colvec{
        \renderStyle \pd{S_{\text{red}}}{N^{(k)}_1} \\\\
        \vdots \\\\
        \renderStyle \pd{S_{\text{red}}}{N^{(n)}_n} \\\\
        \renderStyle \pd{S_{\text{red}}}{V^{(k)}}  
        }  = \colvec{
                   \renderStyle -\frac{\partial \mu_1(T, V^{(k)}, \NBold{k})}{\partial T} + \frac{\partial \mu_1}{\partial T}(\xNew^\Res) \\\\ 
                    \vdots \\\\
                   \renderStyle  -\frac{\partial \mu_n(T, V^{(k)}, \NBold{k})}{\partial T} +\frac{\partial \mu_n}{\partial T}(\xNew^\Res) \\\\
                  \renderStyle  \frac{\partial P(T, V^{(k)}, \NBold{k})}{\partial T} -\frac{\partial P}{\partial T}(\xNew^\Res) 
                    }.
\end{align}
It is interesting to note that, while Smejkal's approach expresses the gradient of \( S_{\text{red}} \) using the terms of the form \( \frac{a}{b} \), our formulation instead involves partial derivatives of the form \( \frac{\partial a}{\partial b} \). For instance, in our approach, the derivative of \(S_{\text{red}}\) with respect to \(N_1^{(k)}\) is expressed as

\[
\pd{S_{\text{red}}}{N_1^{(k)}} = -\pd{\mu_1(T, V^{(k)}, \NBold{k})}{T} + \pd{\mu_1(\xNew^\Res)}{T},
\]
whereas in Smejkal's approach, it is given by

\[
\pd{S_{\text{red}}}{N_1^{(k)}} = -\frac{\mu_1(U^{(k)}, V^{(k)}, \NBold{k})}{T} + \frac{\mu_1(\xNew^\Res)}{T}.
\]
This pattern persists across other derivatives as well, underscoring a fundamental difference in the treatment of thermodynamic variable dependencies between the two methodologies. We remark here that \(U^{(k)}\) is defined differently for both the approaches. For Smejkal's approach, it is the unknown of the optimization problem whereas for our approach it is defined as \(U^{(k)} := U(T, V^{(k)}, \NBold{k})\). Now, having computed the gradients of the entropy function, we now turn our attention to the computation of the gradient of the constraint function \(\cons\).
\begin{align} \label{eqn:cons_grad}    
\nabla \cons({\x}) = \begin{pmatrix}
\renderStyle \nabla \cons^{(1)} \\ 
\vdots \\ 
\renderStyle \nabla \cons^{(p-2)} \\
\renderStyle \nabla \cons^{(p-1)} \\
\renderStyle \pd{\cons}{T}
\end{pmatrix},
\end{align}
\(\forall k \in \{1, \dots, p-1\}, \nabla \cons^{(k)} \in \mathbb{R}^{n+1}\) and \(\pd{\cons}{T} = \sum_{k=1}^{p} \pd{U^{(k)}}{T} \in \mathbb{R}, U^{(k)} := U(T, V^{(k)}, \NBold{k})\).
\begin{align}
    \nabla \cons^{(k)} = \colvec{
        \vertVec{\pd{\cons}{N^{(k)}_1}}{\pd{\cons}{N^{(k)}_n}}{\\\\} \\\\
        \renderStyle \pd{\cons}{V^{(k)}} 
        } = \colvec{
        \vertVec{\pd{\cons^{(k)}}{N^{(k)}_1} - \pd{\cons^{\Res}}{N^{(k)}_1}}{\pd{\cons^{(k)}}{N^{(k)}_n} - \pd{\cons^{\Res}}{N^{(k)}_n}}{\\\\} \\\\
        \renderStyle \pd{\cons^{(k)}}{V^{(k)}} - \pd{\cons^{\Res}}{V^{(k)}} 
        },
\end{align}
We can simplify these gradients further using standard thermodynamic identities. First, recall that the heat capacity at constant volume, \( C_v \), is given by the following thermodynamic relation:
\begin{align} 
    C_v = \left(\pd{U}{T}\right)_{V, \mathbf{N}}.
\end{align}
Consequently, the partial derivative of the constraint with respect to temperature becomes:
\begin{align} \label{eqn:dCdT}
    \pd{\cons}{T} = \sum_{k=1}^{p} C_v^{(k)}.
\end{align}
Next, utilizing the thermodynamic identity:
\begin{align}
    \left( \pd{U}{V} \right)_{T, \mathbf{N}} = T \left( \pd{P}{T} \right)_{V, \mathbf{N}} - P.
\end{align}
we obtain the following expression for the partial derivative of the constraint with respect to volume:
\begin{align}
    \pd{\cons^{(k)}}{V^{(k)}} = \left(\pd{U^{(k)}}{V^{(k)}}\right)_{T, \mathbf{N}} = T \left( \pd{P}{T} \right)_{V^{(k)}, \mathbf{N}} - P.
\end{align}
For the partial derivative of the constraint with respect to the mole number of component 1 in phase \(k\), we proceed as follows:
\begin{align}
\pd{\cons^{(k)}}{N_1^{(k)}} &= \left(\pd{U^{(k)}}{N_1^{(k)}}\right)_{T, V^{(k)}} = \left(\pd{A^{(k)} + T S^{(k)}}{N_1^{(k)}}\right)_{T, V^{(k)}} \nonumber \\
&= \left(\pd{A^{(k)}}{N_1^{(k)}}\right)_{T, V^{(k)}} + T \left(\pd{S^{(k)}}{N_1^{(k)}}\right)_{T, V^{(k)}} \nonumber \\
&= \mu_1^{(k)} - T \left(\pd{\mu_1^{(k)} }{T}\right) ,
\end{align}
where \( \mu_1^{(k)} \) is the chemical potential of component 1 in phase $k$. Similar expressions hold for components \(2, \dots, p-1\).
Finally, we can summarize the gradient of the constraint with respect to the generalized state variables:

\begin{align} \label{eqn:cons_grad}
\nabla \cons^{(k)}({\x}) =
\begin{pmatrix}
\vertVec{\mu^{(k)}_1 - T \pd{\mu^{(k)}_1}{T} - \left( \mu_1^{\Res} - T \pd{\mu_1^{\Res}}{T} \right)}{\mu^{(k)}_n - T \pd{\mu^{(k)}_n}{T} - \left( \mu_n^{\Res} - T \pd{\mu_n^{\Res}}{T} \right)}{\\\\} 
\\\\
\displaystyle T \left( \pd{P^{(k)}}{T} \right)_{V^{(k)}, \NBoldNoSuper} - P^{(k)} - \left(T \left( \pd{P^{\Res}}{T} \right)_{V^{\Res}, \NBoldNoSuper} - P^{\Res} \right) 
\end{pmatrix}.
\end{align}
With the gradient of the entropy and constraint  fully defined, we now proceed to the computation of the Lagrange multiplier. This multiplier plays a crucial role in enforcing the constraint and is directly determined by the relationship between these gradients.

\subsection{Hessian computation}
The Hessian matrix \( H(\xNew) \in \mathbb{R}^{(p-1)(n+1) + 1, (p-1)(n+1) + 1} \) is given in block form as:
\[
\HessNew(\xNew) =
\begin{bmatrix}
    \HessNew_{N,N} & \HessNew_{N,V} & \HessNew_{N,T} \\
    \HessNew_{V,N} & \HessNew_{V,V} & \HessNew_{V,T} \\
    \HessNew_{T,N} & \HessNew_{T,V} & \HessNew_{T,T}
\end{bmatrix}
\]

Each block has the following structure:

\(\HessNew_{N,N}\) (\((p-1)n \times (p-1)n\)):
\[
\HessNew_{N,N}^{(k,\ell)} =
\begin{bmatrix}
    \renderStyle \frac{\partial^2 S}{\partial N_1^{(k)} \partial N_1^{(\ell)}} & \cdots & \renderStyle \frac{\partial^2 S}{\partial N_1^{(k)} \partial N_n^{(\ell)}} \\
    \renderStyle \vdots & \ddots & \vdots \\
    \renderStyle \frac{\partial^2 S}{\partial N_n^{(k)} \partial N_1^{(\ell)}} & \cdots & \renderStyle \frac{\partial^2 S}{\partial N_n^{(k)} \partial N_n^{(\ell)}}
\end{bmatrix}
\]

\(\HessNew_{N,V}\) (\((p-1)n \times (p-1)\)) and \(\HessNew_{N,T}\) (\((p-1)n \times 1\)):
\[
\HessNew_{N,V}^{(k,\ell)} =
\begin{bmatrix}
    \renderStyle \frac{\partial^2 S}{\partial N_1^{(k)} \partial V^{(\ell)}} \\
    \vdots \\
    \renderStyle \frac{\partial^2 S}{\partial N_n^{(k)} \partial V^{(\ell)}}
\end{bmatrix}, \quad
\HessNew_{N,T}^{(k)} =
\begin{bmatrix}
    \renderStyle \frac{\partial^2 S}{\partial N_1^{(k)} \partial T} \\
    \vdots \\
    \renderStyle \frac{\partial^2 S}{\partial N_n^{(k)} \partial T}
\end{bmatrix}
\]

\(\HessNew_{V,V}\) (\((p-1) \times (p-1)\)) and \(\HessNew_{V,T}\) (\((p-1) \times 1\)):
\[
\HessNew_{V,V}^{(k,\ell)} = \frac{\partial^2 S}{\partial V^{(k)} \partial V^{(\ell)}}, \quad
\HessNew_{V,T}^{(k)} = \frac{\partial^2 S}{\partial V^{(k)} \partial T}
\]

\(\HessNew_{T,T}\) (\(1 \times 1\)):
\[
\HessNew_{T,T} = \frac{\partial^2 S}{\partial T^2}
\]
Similarly, the Hessian of the constraint function can be derived using the same methodology. In practice, however, we leverage state-of-the-art automatic differentiation (AD) tools to compute these derivatives accurately. While AD provides a robust and reliable framework for evaluating gradients and Hessians, manual derivations remain invaluable for verifying the correctness of AD-generated results. Additionally, although AD is generally efficient, computing higher-order derivatives can become computationally intensive. In such scenarios, manually derived expressions for derivatives may offer a more computationally efficient alternative. In this work, we primarily rely on the \texttt{ForwardDiff.jl} package for AD to compute the gradient and the Hessian, and use manual derivations for verification purposes.

\section{Lagrange Multiplier} \label{app:lagrange_mul}

The \textbf{Lagrange multiplier} method is used to find the extrema of a function subject to constraints. Given an objective function \( f(x_1, x_2, \dots, x_n) \) and a constraint 

\[
g(x_1, x_2, \dots, x_n) = 0,
\]

we define the \textbf{Lagrangian} function:

\[
\mathcal{L}(x_1, x_2, \dots, x_n, \lambda) = f(x_1, x_2, \dots, x_n) + \lambda g(x_1, x_2, \dots, x_n).
\]

The necessary conditions for an extremum are obtained by setting the gradient of \( \mathcal{L} \) to zero:

\[
\nabla \mathcal{L} = 0 \quad \Rightarrow \quad \nabla f = -\lambda \nabla g.
\]

Additionally, the constraint equation must be satisfied:

\[
g(x_1, x_2, \dots, x_n) = 0.
\]

These conditions ensure that the gradient of \( f \) is parallel to the gradient of \( g \), meaning the constraint surface is tangent to the level curves of \( f \).

\section{Peng-Robinson Equation of State} \label{app:PR_EOS}

We employ the Peng-Robinson equation of state (EOS) \cite{smejkal_phase_2017}, which is formulated as follows:
\begin{equation}
    P^{(\text{EOS})}(T, V, N_1, \ldots, N_n) = \frac{NRT}{V - Nb} - \frac{a(T)N^2}{V^2 + 2bNV - N^2b^2}, \label{eq:peng_robinson}
\end{equation}
where \(T\) is the temperature, \(V\) is the volume, \(N_i\) represents the number of moles of component \(i\) in the system, \(R\) is the universal gas constant and \(N\) is the total number of moles in the system. The parameters \(a(T)\) and \(b\) characterize intermolecular forces and volume exclusion, respectively. The parameters \(a(T)\) and \(b\) are defined as follows:
\begin{subequations}
\begin{equation}
    a = \sum_{i=1}^{n} \sum_{j=1}^{n} x_i x_j a_{ij}, \label{eq:a_parameter}
\end{equation}
\begin{equation}
    a_{ij} = (1 - \delta_{ij}) \sqrt{a_i a_j}, \label{eq:a_ij}
\end{equation}
\begin{equation}
    a_i(T) = 0.45724 \frac{R^2 T_{\text{crit},i}^2}{P_{\text{crit},i}} \left[1 + m_i \left(1 - \sqrt{T_{r,i}} \right)\right]^2, \label{eq:a_i}
\end{equation}
\begin{equation}
    b = \sum_{i=1}^{n} x_i b_i, \label{eq:b_parameter}
\end{equation}
\begin{equation}
    b_i = 0.0778 \frac{R T_{\text{crit},i}}{P_{\text{crit},i}}, \label{eq:b_i}
\end{equation}
\end{subequations}
where \(x_i\) is the mole fraction of component \(i\), \(T_{\text{crit},i}\) and \(P_{\text{crit},i}\) are the critical temperature and pressure of component \(i\), and \(\delta_{ij}\) is the Kronecker delta. The parameter \(m_i\) accounts for the acentric factor \(\omega_i\) as:
\begin{equation}
    m_i = \begin{cases} 
    0.37464 + 1.54226 \omega_i - 0.26992 \omega_i^2, & \omega_i < 0.5, \\
    0.3796 + 1.485 \omega_i - 0.1644 \omega_i^2 + 0.01667 \omega_i^3, & \omega_i \geq 0.5.
    \end{cases} \label{eq:m_i}
\end{equation}

The internal energy, \(U^{(\text{EOS})}\), and entropy, \(S^{(\text{EOS})}\), in the context of the Peng-Robinson EOS are expressed as follows:

The internal energy:
\begin{align}
    U^{(\text{EOS})}(T, V, N_1, \ldots, N_n) &= N \frac{T \partial_T(a) - a} {2 \sqrt{2}b} \ln \left| \frac{V + (1 + \sqrt{2})bN}{V + (1 - \sqrt{2})bN} \right| \nonumber \\
    &- NR(T - T_0) + \sum_{i=1}^{n} N_i \sum_{k=0}^{3} \alpha_{ik} \frac{T^{k+1} - T_0^{k+1}}{k+1} + N u_0, \label{eq:internal_energy}
\end{align}
where \(\partial_T(a)\) is the temperature derivative of \(a(T)\), \(T_0\) is a reference temperature, and \(\alpha_{ik}\) are empirical constants.

The entropy:
\begin{align}
    S^{(\text{EOS})}(T, V, N_1, \ldots, N_n) &= NR \ln \left| \frac{V - bN}{V} \right| + N\frac{\partial_T(a)}{2 \sqrt{2}b} \ln \left| \frac{V + (1 + \sqrt{2})bN}{V + (1 - \sqrt{2})bN} \right| \nonumber \\
    &+ R \sum_{i=1}^{n} N_i \ln \frac{V P_0}{N_i RT} + \sum_{i=1}^{n} N_i \int_{T_0}^{T} \frac{c_{p,i}^{\text{ig}}(\xi)}{\xi} d\xi, \label{eq:entropy}
\end{align}
where \(c_{p,i}^{\text{ig}}(T)\) is the ideal gas heat capacity of component \(i\) and \(P_0\) is a reference pressure. The heat capacity \(c_{p,i}^{\text{ig}}(T)\) can be written as:
\begin{equation}
    c_{p,i}^{\text{ig}}(T) = \sum_{k=0}^{3} \alpha_{ik} T^k. \label{eq:heat_capacity}
\end{equation}

\end{document}